\documentclass[11pt]{article}

\usepackage{amsmath}
\usepackage{amssymb}
\usepackage[mathcal]{euscript}
\usepackage{graphicx,psfrag}
\usepackage{epsfig}
\usepackage{epic,eepic}

\newcommand{\ee}{{\bf e}}

\newcommand{\bbZ}{{\mathbb Z}}
\newcommand{\bbR}{{\mathbb R}}

\newcommand{\bbC}{{\mathbb C}}
\newcommand{\bbL}{{\mathbb L}}
\newcommand{\bbS}{{\mathbb S}}
\newcommand{\bbP}{{\mathbb P}}

\newcommand{\cC}{{\mathcal C}}

\newcommand{\cL}{{\mathcal L}}

\newcommand{\cQ}{{\mathcal Q}}

\newcommand{\cZ}{{\mathcal Z}}

\newtheorem{theorem}{Theorem}
\newtheorem{definition}[theorem]{Definition}
\newtheorem{corollary}[theorem]{Corollary}
\newtheorem{lemma}[theorem]{Lemma}

\begin{document}
\title{On Organizing Principles of
\\
Discrete Differential Geometry. \\
Geometry of spheres}

\author{
Alexander I. Bobenko
\thanks{Institut f\"ur Mathematik,
Technische Universit\"at Berlin, Str. des 17. Juni 136, 10623
Berlin, Germany. E--mail: {\tt bobenko@math.tu-berlin.de}} \and
Yuri B. Suris
\thanks{Zentrum Mathematik, Technische Universit\"at M\"unchen,
Boltzmannstr. 3, 85747 Garching bei M\"unchen, Germany. E--mail:
{\tt suris@ma.tum.de}}}

\maketitle{\renewcommand{\thefootnote}{} \footnote[0]{Research for
this article was supported by the DFG Research Unit 565
`Polyhedral Surfaces'' and the DFG Research Center
\textsc{Matheon} ``Mathematics for key technologies'' in Berlin.}}

{\small {\bf Abstract.}
  Discrete differential geometry aims to develop discrete equivalents
of the geometric notions and methods of classical differential
geometry. In this survey we discuss the following two fundamental
Discretization Principles: the transformation group principle
(smooth geometric objects and their discretizations are invariant
with respect to the same transformation group) and the consistency
principle (discretizations of smooth parametrized geometries can
be extended to multidimensional consistent nets). The main
concrete geometric problem discussed in this survey is a
discretization of curvature line parametrized surfaces in Lie
geometry. By systematically applying the Discretization Principles
we find a discretization of curvature line parametrization which
unifies the circular and conical nets.}

\newpage
\section{Introduction}
\label{Sect: intro}

The new field of discrete differential geometry is presently
emerging on the border between differential and discrete geometry.
Whereas classical differential geometry investigates smooth
geometric shapes, discrete differential geometry studies geometric
shapes with finite numbers of elements and aims to develop
discrete equivalents of the geometric notions and methods of
classical differential geometry. The latter appears then as a
limit of refinements of the discretization. Current interest in
this field derives not only from its importance in pure
mathematics but also from its relevance for computer graphics. An
important example one should keep in mind here are polyhedral
surfaces approximating smooth surfaces.


One may suggest many different reasonable discretizations with the
same smooth limit. Which one is the best? From the theoretical
point of view the best discretization is the one which preserves
all fundamental properties of the smooth theory. Often such a
discretization clarifies the structures of the smooth theory and
possesses important connections to other fields of mathematics
(projective geometry, integrable systems, algebraic geometry,
complex analysis etc.). On the other hand, for applications the
crucial point is the approximation: the best discretization is
supposed to possess distinguished convergence properties and
should represent a smooth shape by a discrete shape with just few
elements. Although these theoretical and applied criteria for the
best discretization are completely different,
in many cases natural ``theoretical'' discretizations turn out to
possess remarkable approximation properties
and are very useful for applications \cite{BSch, LPWYW}.


This interaction of the discrete and smooth versions of the theory
led to important results in the surface theory as well as in the
geometry of polyhedra. Classical achievements of discrete
differential geometry are the fundamental results of Alexandrov
and Pogorelov on metric geometry of polyhedra and convex surfaces:
Alexandrov's theorem \cite{Alexandrov} states that any abstract
convex polyhedral metric is uniquely realized by a convex
polyhedron in Euclidean 3-space. Pogorelov proved \cite{Pogorelov}
the corresponding existence and uniqueness result for abstract
convex metrics by approximating smooth surfaces by polyhedra.


Simplicial surfaces, i.e., discrete surfaces made from triangles,
are basic in computer graphics. This class of discrete surfaces,
however, is too unstructured for analytical investigation. An
important tool in the theory of smooth surfaces is the
introduction of (special) parametrizations of a surface. Natural
analogues of parametrized surfaces are {\em quadrilateral
surfaces}, i.e. discrete surfaces made from (not necessarily
planar) quadrilaterals. The strips of quadrilaterals obtained by
gluing quadrilaterals along opposite edges are analogs of
coordinate lines. Probably the first nontrivial example of
quadrilateral surfaces studied this way are discrete surfaces with
constant negative Gaussian curvature introduced by Sauer and
Wunderlich \cite{SauerK, Wu}. Currently discrete parametrized
surfaces are becoming more important in computer graphics. They
lead to meshes that better represent the shape of the surface and
look regular \cite{ACSDLD, DKG, MK, LPWYW}.


It is well known that differential equations describing
interesting special classes of surfaces and parametrizations are
integrable (in the sense of the theory of integrable systems),
and, conversely, many of interesting integrable systems admit a
differential-geometric interpretation. A progress in understanding
of the unifying fundamental structure the classical differential
geometers were looking for, and simultaneously in understanding of
the very nature of integrability, came from the efforts to
discretize these theories. It turns out that many sophisticated
properties of differential-geometric objects find their simple
explanation within the discrete differential geometry. The early
period of this development is documented in the work of Sauer
\cite{Sauer}. The modern period began with the work by Bobenko and
Pinkall \cite{BP1, BP2} and by Doliwa and Santini \cite{DS1, CDS}.
A closely related development of the spectral theory of difference
operators on graphs was initiated by Novikov with collaborators
\cite{ND, No1, No2}, see also \cite{DN} for a further development
of a discrete complex analysis on simplicial manifolds.

Discrete surfaces in Euclidean 3-space is the basic example
considered in this survey. This case has all essential features of
the theory in all generality, generalizations for higher
dimensions are straightforward. On the other hand, our geometric
three-dimensional intuition helps to understand their properties.

Discrete differential geometry related to integrable systems deals
with multidimensional discrete nets, i.e., maps from the regular
cubic lattice ${\mathbb Z}^m$ into ${\mathbb R}^N$ specified by
certain geometric properties (as mentioned above, we will be most
interested in the case $N=3$ in this survey). In this setting
discrete surfaces appear as two dimensional layers of
multidimensional discrete nets, and their transformations
correspond to shifts in the transversal lattice directions. A
characteristic feature of the theory is that all lattice
directions are on equal footing with respect to the defining
geometric properties. Discrete surfaces and their transformations
become indistinguishable. We associate such a situation with the
{\em multidimensional consistency}, and this is one of our
fundamental discretization principles. The multidimensional
consistency, and therefore the existence and construction of
multidimensional nets, relies just on certain incidence theorems
of elementary geometry.

\begin{figure}[htbp]
\begin{center}
\includegraphics[width=0.40\textwidth]{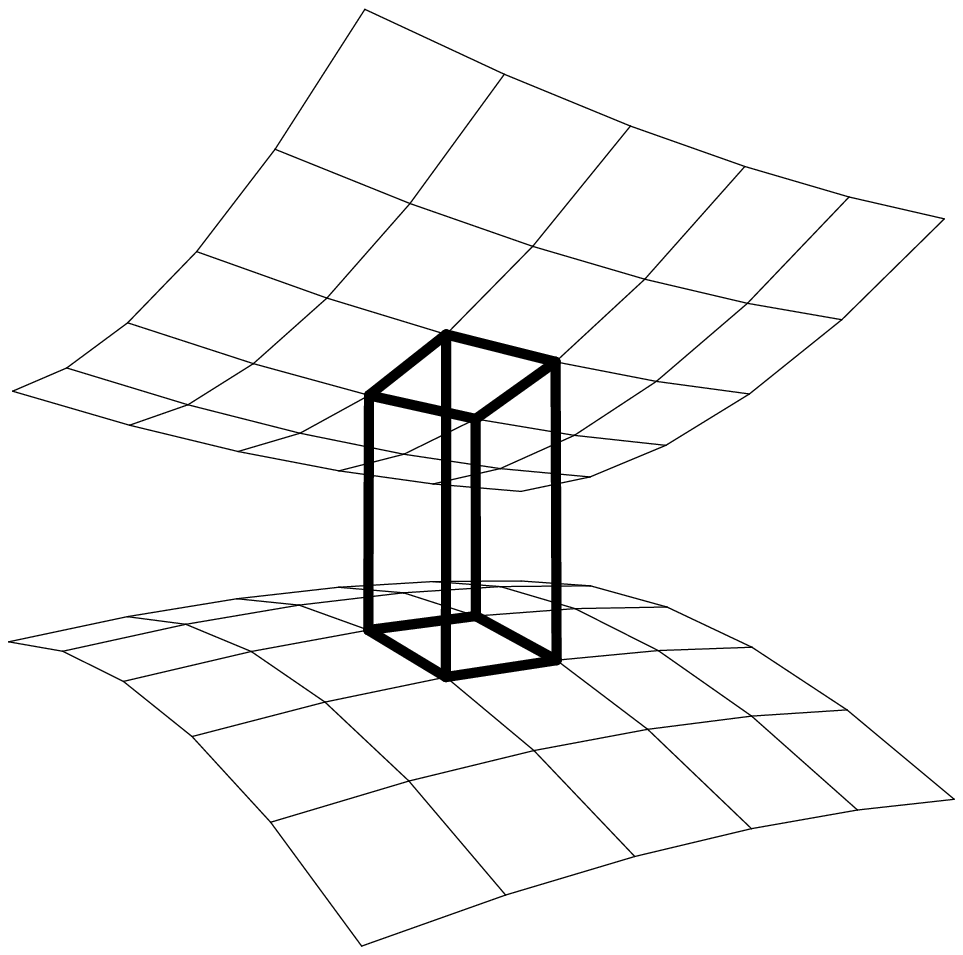}\qquad
\includegraphics[width=0.40\textwidth]{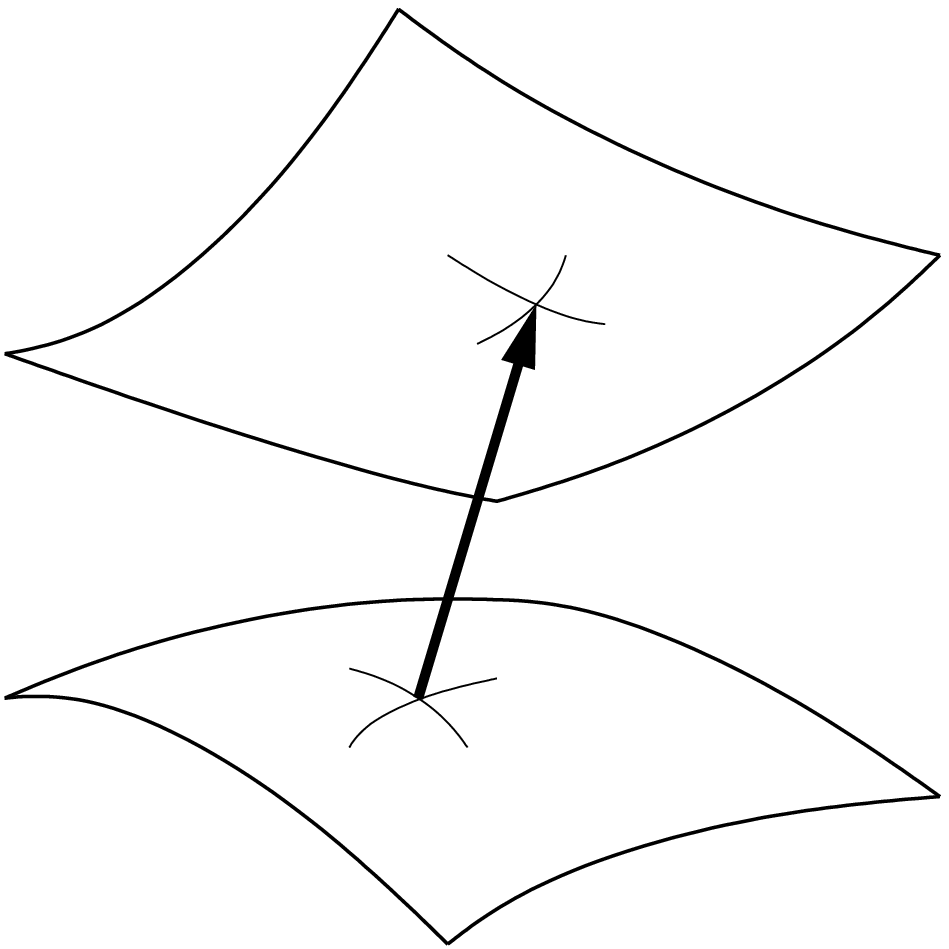}
\end{center}
\caption{From the discrete master theory to the classical theory:
surfaces and their transformations appear by refining two of three
net directions.} \label{fig:discrete->classical}
\end{figure}

Conceptually one can think of passing to a continuum limit by
refining mesh size in some of the lattice directions. In these
directions the net converges to smooth surfaces whereas those
directions that remain discrete correspond to transformations of
the surfaces (see Fig. \ref{fig:discrete->classical}).
The smooth theory comes as a corollary of a more fundamental
discrete master theory. The true roots of the classical surface
theory are found, quite unexpectedly, in various incidence
theorems of elementary geometry. This phenomenon, which has been
shown for many classes of surfaces and coordinate systems
\cite{BP2, BSDDG}, currently is getting accepted as one of the
fundamental features of classical integrable differential
geometry.

Note that finding simple discrete explanations for complicated
differential geometric theories is not the only outcome of this
development. Having identified the roots of the integrable
differential geometry in the multidimensional consistency of
discrete nets, we are led to a new (geometric) understanding of
the integrability itself \cite{BS1, ABS, BSDDG}.


The simplest and at the same time the basic example of consistent
multidimensional nets are multidimensional Q-nets \cite{DS1}, or
discrete conjugate nets \cite{Sauer}, which are characterized by
planarity of all quadrilaterals. The planarity property is
preserved by projective transformations and thus Q-nets are
subject of projective geometry (like conjugate nets, which are
smooth counterparts of Q-nets).

Here we come to the next basic discretization principle. According
to F. Klein's Erlangen program, geometries are classified by their
transformation groups. Classical examples are projective, affine,
Euclidean, spherical, hyperbolic geometry, and sphere geometries
of Lie, M\"oebius, and Laguerre. We postulate that the
transformation group as the most fundamental feature should be
preserved by a discretization. This can be seen as a sort of {\em
discrete Erlangen program}.

Thus we come to the following fundamental \smallskip

{\bf Discretization Principles}:

\begin{itemize}
 \item {\em Transformation group principle:} smooth geometric objects and
    their discretizations belong to the same geometry, i.e.
    are invariant with respect to the same transformation group.
 \item {\em Consistency principle:} discretizations of smooth
    parametrized geometries can be extended to multidimensional
    consistent nets.
\end{itemize}


Let us explain why such different imperatives as the
transformation group principle and the consistency principle can
be simultaneously imposed for discretization of classical
geometries. The transformation groups of various geometries,
including those of Lie, M\"obius and Laguerre, are subgroups of
the projective transformation group. Classically, such a subgroup
is described as consisting of projective transformations which
preserve some distinguished quadric called absolute. A remarkable
result by Doliwa \cite{D2} is that multidimensional Q-nets can be
restricted to an arbitrary quadric. This is the reason why the
Discretization Principles work for the classical geometries.

In this survey we deal with three classical geometries described
in terms of spheres: M\"obius, Laguerre and Lie geometries. They
have been developed by classics, the most elaborated presentation
of these geometries can be found in Blaschke's book \cite{BlIII}.

M\"obius geometry is the most popular one of these three
geometries. It describes properties invariant with respect to
M\"obius transformations which are compositions of reflections in
spheres. For $N\ge 3$, M\"obius transformations of $\bbR^N$
coincide with conformal transformations. M\"obius geometry does
not distinguish between spheres and planes (planes are regarded as
spheres through the infinitely remote point $\infty$, which
compactifies $\bbR^N$ to the $N$-sphere $\bbS^N$). On the other
hand, points are considered as objects different from spheres.
Surfaces are described through their points. Classical examples of
M\"obius-geometric properties of surfaces are conformal
parametrization and the Willmore functional \cite{Willmore}.
Recent progress in this field is to a large extent due to
interrelations with the theory of integrable systems \cite{FLPP,
Taimanov}.

Laguerre geometry does not distinguish points and spheres (points
are treated as spheres of zero radius). On the other hand, planes
are distinguished. Surfaces are described through their tangent
planes. A particular Laguerre transformation of a surface is a
shift of all tangent planes in the normal direction at a constant
distance. This transformation is called the normal shift.

Lie geometry is a natural unification of M\"obius and Laguerre
geometries: points, planes and spheres are treated on an equal
footing. The transformation group is generated by M\"obius
transformations and the normal shift transformations. Surfaces are
described through their contact elements. A contact element can be
understood as a  surface point together with the corresponding
tangent plane. The one-parameter family of spheres through a point
with a common tangent plane gives a Lie-geometric description of a
contact element. The point of the surface and the tangent plane at
this point are just two elements of this family.

Integrability aspects of the surface theory in Lie geometry have
been studied by Ferapontov \cite{F1, F2}, Musso and Nicolodi
\cite{MN}, and Burstall and Hertrich-Jeromin \cite{BHJ1, BHJ2}.


The main concrete geometric problem discussed in this survey is a
discretization of curvature line parametrized surfaces. Curvature
lines are integral curves of the principal directions. Any surface
away from its umbilic points can be parametrized by curvature
lines. Curvature line parametrization is attracting the attention
of mathematicians and physicists for two centuries. The classical
results in this field can be found in the books by Darboux
\cite{Da1, Da2} and Bianchi \cite{Bi}. In particular a classical
result of Dupin \cite{Dupin} claims that the coordinate surfaces
of triply orthogonal coordinate systems intersect along their
common curvature lines. Ribaucour has discovered a transformation
of surfaces preserving the curvature line parametrization (see
\cite{E2}). A surface and its Ribaucour transform envelope a
special sphere congruence. Bianchi has shown \cite{Bi2} that
Ribaucour transformations are permutable: given two Ribaucour
transforms of a surface there exists a one-parameter family of
their common Ribaucour transforms.

Recently curvature line parametrizations and orthogonal systems
came back into the focus of interest in mathematical physics as an
example of an integrable system. Zakharov \cite{Z} has constructed
a variety of explicit solutions with the help of the dressing
method. Algebro-geometric orthogonal coordinate systems were
constructed by Krichever \cite{Kri}. The recent interest to this
problem is in particular motivated by applications to the theory
of the associativity equations developed by Dubrovin \cite{Dub}.
Remarkable geometric properties make curvature line
parametrizations especially useful for visualization of surfaces
in computer graphics \cite{ACSDLD, LPWYW}.

The question of proper discretization of the curvature line
parametrized surfaces and orthogonal systems became recently a
subject of intensive study. {\em Circular nets}, which are Q-nets
with circular quadrilaterals, as discrete analogs of curvature
line parametrized surfaces were mentioned by Nutbourne and Martin
\cite{NM}. Special circular nets as discrete isothermic surfaces
were investigated in \cite{BP2}. The circular discretization of
triply-orthogonal coordinate systems was first suggested by one of
the authors in \cite{B}. Doliwa and Santini \cite{DS1} made the
next crucial step in the development of the theory. They
considered discrete orthogonal systems as a reduction of discrete
conjugated systems \cite{CDS}, generalized them to arbitrary
dimension and proved their multidimensional consistency based on
the classical Miquel theorem \cite{Be}.

Matthes and the authors of this survey have proven \cite{BMS2}
that circular nets approximate smooth curvature line parametrized
surfaces and orthogonal systems with all derivatives. Numerical
experiments show that circular nets have the desired geometrical
properties already at the coarse level and not only in the
refinement limit as it approaches a smooth curvature line
parameterized surface. This is important for applications in
computer graphics \cite{LPWYW}.

A convenient analytic description of circular nets has been given
by Konopelchenko and Schief \cite{KSch}. Analytic methods of the
soliton theory have been applied to circular nets by Doliwa,
Manakov and Santini \cite{DMS} ($\bar\partial$-method) and by
Akhmetshin, Volvovskii and Krichever \cite{AKV} (algebro-geometric
solutions). Bobenko and Hertrich-Jeromin \cite{BJ} have given a
Clifford algebra description of circular nets.

Circular nets are preserved by M\"obius transformations, and thus
should be treated as the discretization of curvature line
parametrizations in M\"obius geometry. A remarkable recent
development by Liu, Pottmann, Wallner, Yang, and Wang \cite{LPWYW}
is the introduction of {\em conical nets}, which should be treated
as the discretization of curvature line parametrizations in
Laguerre geometry. These are special Q-nets characterized by the
property that four quadrilaterals meeting at a vertex are tangent
to a common cone of revolution. Equivalently, conical nets can be
characterized as Q-nets with circular Gauss maps, i.e., the unit
normals to the quadrilaterals comprise a circular net in the unit
sphere $S^2$. Circular Gauss maps defined at vertices of a given
circular net were previously introduced by Schief \cite{Sch2},
however without relation to conical nets. Conical nets, like
circular ones, satisfy the second discretization principle
(consistency).

In the present survey, we find a discretization of curvature line
parametrization which unifies the circular and the conical nets by
systematically applying the Discretization Principles.

\begin{figure}[hb]
\begin{center}
\begin{picture}(0,0)%
\includegraphics{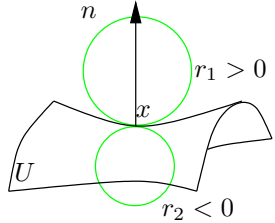}%
\end{picture}%
\setlength{\unitlength}{4144sp}%
\begingroup\makeatletter\ifx\SetFigFont\undefined%
\gdef\SetFigFont#1#2#3#4#5{%
  \reset@font\fontsize{#1}{#2pt}%
  \fontfamily{#3}\fontseries{#4}\fontshape{#5}%
  \selectfont}%
\fi\endgroup%
\begin{picture}(1599,1320)(259,-471)
\put(706,714){\makebox(0,0)[lb]{\smash{\SetFigFont{10}{12.0}{\rmdefault}{\mddefault}
{\updefault}
$n$
}}}
\put(1191,-471){\makebox(0,0)[lb]{\smash{\SetFigFont{10}{12.0}{\rmdefault}{\mddefault}
{\updefault}
$r_2<0$
}}}
\put(311,-261){\makebox(0,0)[lb]{\smash{\SetFigFont{10}{12.0}{\rmdefault}{\mddefault}
{\updefault}
$U$
}}}
\put(1386,364){\makebox(0,0)[lb]{\smash{\SetFigFont{10}{12.0}{\rmdefault}{\mddefault}
{\updefault}
$r_1>0$
}}}
\put(1036,119){\makebox(0,0)[lb]{\smash{\SetFigFont{10}{12.0}{\rmdefault}{\mddefault}
{\updefault}
$x$
}}}
\end{picture}
\end{center}
\caption{Principal directions through touching spheres.}
\label{fig:p_directions}
\end{figure}

It is well known that curvature lines are subject of Lie geometry,
i.e., are invariant with respect to M\"obius transformations and
normal shifts. To see this, consider an infinitesimal neighborhood
$U$ of a point $x$ of an oriented smooth surface in $\bbR^3$, and
the pencil of spheres $S(r)$ of the signed radii $r$, touching the
surface at $x$, see Fig.~\ref{fig:p_directions}. The signed radius
$r$ is assumed positive if $S(r)$ lies on the same side of the
surface as the normal $n$, and negative otherwise; $S(\infty)$ is
the tangent plane. For small $r_0>0$ the spheres $S(r_0)$ and
$S(-r_0)$ intersect $U$ in $x$ only. The set of the touching
spheres with this property (intersecting $U$ in $x$ only) has two
connected components: $M_+$ containing $S(r_0)$ and $M_-$
containing $S(-r_0)$ for small $r_0>0$. The boundary values
\begin{eqnarray*}
r_1=\sup\{r: S(r)\in M_+\},\qquad r_2=\inf\{r: S(r)\in M_-\}
\end{eqnarray*}
are the principal curvatures of the surface in $x$. The directions
in which $S(r_1)$ and $S(r_2)$ touch $U$ are the principal
directions.

Clearly, all ingredients of this description are
M\"obius-invariant. Under a normal shift by the distance $d$ the
centers of the principal curvature spheres are preserved and their
radii are shifted by $d$. This implies that the principal
directions and thus the curvature lines are preserved under normal
shifts, as well.

A Lie-geometric nature of the curvature line parametrization
yields that it has a Lie-invariant description. Such a description
can be found in Blaschke's book \cite{BlIII}. A surface in Lie
geometry, as already said, is considered as consisting of contact
elements. Two infinitesimally close contact elements (sphere
pencils) belong to the same curvature line, if and only if they
have a sphere in common, which is the principal curvature sphere.

By a literal discretization of this Blaschke's Lie-geometric
description of smooth curvature line parametrized surfaces, we
define a discrete {\em principal contact element net} as a map
$\bbZ^2\to\{{\rm contact\ elements\ of\ surfaces\ in\ }\bbR^3\}$
such that any two neighboring contact elements have a sphere in
common.

In the projective model of Lie geometry spheres in $\bbR^3$
(including points and planes) are represented by elements of the
so called Lie quadric $\bbL\subset\bbR\bbP^{5}$, contact elements
are represented by isotropic lines, i.e., lines in $\bbL$,
surfaces are represented by congruences of isotropic lines. In the
curvature line parametrization, the parametric families of
isotropic lines comprise developable surfaces in $\bbL$.

Accordingly, a discrete principal contact element net in the
projective model of Lie geometry is a discrete congruence of
isotropic lines
\[
\ell:\bbZ^2\to\{\rm isotropic\ lines\ in\ \bbL\}
\]
such that any two neighboring lines intersect. Intersection points
of neighboring lines are, as in the smooth case, the principal
curvature spheres. They are associated with the edges of $\bbZ^2$.
Four principal curvature spheres associated to the edges with a
common vertex belong to the same contact element, i.e., have a
common touching point.
\begin{figure}[htbp]
\begin{center}
\rotatebox{-90}{\includegraphics[width=0.5\textwidth]{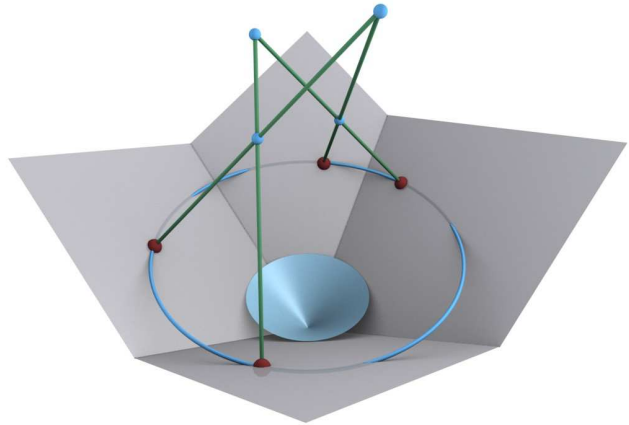}}
\end{center}
\caption{Geometry of principal contact element nets. Four
neighboring contact elements are represented by points and
(tangent) planes. The points are concircular, the planes are
tangent to a cone of revolution. Neighboring normal lines
intersect at the centers of principal curvature spheres.}
\label{fig:Lie-curvature}
\end{figure}

In projective geometry, discrete line congruences have been
introduced by Doliwa, Ma\~nas and Santini \cite{DMS}. Discrete
line congruences are closely related to Q-nets, and, like the
latter, are multidimensionally consistent. It follows from our
results that they can be restricted to the Lie quadric (actually,
to any ruled quadric). Thus, principal contact element nets
satisfy the second discretization principle. In particular, this
yields discrete Ribaucour transformations between principal
contact element nets.

The Lie-geometric notion of discrete principal contact element
nets unifies the M\"obius-geometric one (circular nets) and the
Laguerre-geometric one (conical nets). Indeed, any contact element
$\ell$ contains a point $x$ and a plane $P$. It turns out that for
a surface
\[
\ell:\bbZ^2\to\{\rm isotropic\ lines\ in\ \bbL\}=\{\rm contact\
elements\ in\ \bbR^3\},
\]
the points comprise a circular net
\[
x:\bbZ^2\to\bbR^3,
\]
whereas the planes comprise a conical net
\[
P:\bbZ^2\to\{\rm planes\ in\ \bbR^3\}.
\]
The corresponding geometry is depicted on
Fig.~\ref{fig:Lie-curvature}.
\begin{figure}[htbp]
\begin{center}
\setlength{\unitlength}{0.055em}
\begin{picture}(200,220)(0,0)
 \put(0,0){\circle*{10}}    \put(150,0){\circle*{10}}
 \put(0,150){\circle*{10}}  \put(150,150){\circle*{10}}
 \put(50,50){\circle*{10}} \put(50,200){\circle*{10}}
 \put(200,50){\circle*{10}}
 \put(200,200){\circle*{10}}
 \path(0,0)(150,0)       \path(0,0)(0,150)
 \path(150,0)(150,150)   \path(0,150)(150,150)
 \path(0,150)(50,200)    \path(150,150)(200,200)
 \path(50,200)(200,200)
 \path(200,200)(200,50) \path(200,50)(150,0)
 \dashline[+30]{10}(0,0)(50,50)
 \dashline[+30]{10}(50,50)(50,200)
 \dashline[+30]{10}(50,50)(200,50)
 \put(-25,-5){$x$}
 \put(-25,145){$P$}
 \put(-25,70){$\ell$}
 \put(23,50){$x_2$}
 \put(22,200){$P_2$}
 \put(25,125){$\ell_2$}
 \put(160,-5){$x_1$}
 \put(160,140){$P_1$}
 \put(160,70){$\ell_1$}
 \put(215,50){$x_{12}$}
 \put(215,200){$P_{12}$}
 \put(215,125){$\ell_{12}$}
 \put(65,25){M\"OBIUS}
 \put(60,170){LAGUERRE}
 \put(165,100){LIE}
\end{picture}
\caption{Geometry of principal contact element nets. Four
neighboring contact elements produce a hexahedron with vertices in
the Lie quadric $\bbL$ and with planar faces. The bottom
quadrilateral is the intersection of the three-dimensional space
$V={\rm span}(\ell,\ell_1,\ell_2,\ell_{12})$ with the 4-space in
$\bbR\bbP^5$ representing points in $\bbR^3$. The top
quadrilateral is the intersection of $V$ with the 4-space in
$\bbR\bbP^5$ representing planes in $\bbR^3$. Each side
quadrilateral lies in the plane of two intersecting lines
$\ell\subset\bbL$.}\label{fig:merging}
\end{center}
\end{figure}
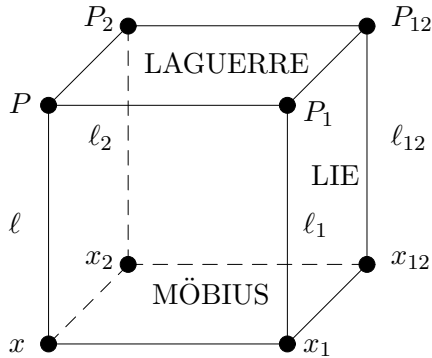
Schematically, this Lie-geometric merging of the M\"obius- and
Laguerre-geometric notions is presented on Fig.~\ref{fig:merging}.

This survey is organized as follows. In Sect.~\ref{Sect: consist}
we start with a review of the basic multidimensionally consistent
systems -- the Q-nets and discrete line congruences. The basic
notions of Lie, M\"obius and Laguerre geometries are briefly
presented in Sect.~\ref{Sect: geom}. Sect.~\ref{Sect: curv}
contains main new results on discrete curvature line parametrized
surfaces: Lie-geometric Definitions \ref{def: Lie curv Euc},
\ref{def: Lie curv}, and Theorem \ref{thm: synthesis} which
describes interrelations of discrete curvature line nets in Lie,
M\"obius and Laguerre geometries. Geometric characterization of
Ribaucour transformations and discrete R-congruences of spheres as
quadrilateral nets in the Lie quadric is given in Sect.~\ref{sect:
R-congr}.

Let us note that, due to the classical Lie's sphere-line
correspondence, the Lie-geometric theory presented in this survey
can be transferred to the context of projective line geometry in
three-space: the Lie quadric is replaced by the Pl\"ucker quadric,
the curvature lines and R-congruences of spheres correspond to the
asymptotic lines and the W-congruences of lines, respectively. The
projective theory of discrete asymptotic nets has been developed
by Doliwa \cite{D4}.

Our research in discrete differential Lie geometry has been
stimulated by the recent introduction of conical nets by Liu,
Pottmann, Wallner, Yang, and Wang \cite{LPWYW}. The advent of the
second (after circular nets) discretization of curvature line
parametrizations posed a question about the relation between the
different discretizations. Independently, a relation between
circular and conical nets has been found by Pottmann \cite{P}. We
are grateful to H.~Pottmann and J.~Wallner for numerous
communications on conical nets and for providing us with their
unpublished results. We thank also U.~Pinkall for useful
discussions.

\section{Consistency as a discretization principle}
\label{Sect: consist}

\subsection{Q-nets}
\label{subsect: Q-nets}

We use the following standard notation: for a function $f$ on
$\bbZ^m$ we write
\[
\tau_i f(u)=f(u+e_i),
\]
where $e_i$ is the unit vector of the $i$-th coordinate direction,
$1\le i\le m$. We use also the shortcut notations $f_i$ for
$\tau_i f$, $f_{ij}$ for $\tau_i\tau_j f$, etc.

The most general of the known discrete 3D systems possessing the
property of 4D consistency are nets consisting of planar
quadrilaterals, or Q-nets. Two-dimensional Q-nets were introduced
by Sauer \cite{Sauer}, the multi-dimensional generalization has
been given by Doliwa and Santini \cite{DS1}. Our presentation in
this section follows the latter paper. The fundamental importance
of multi-dimensional consistency of discrete systems as their
integrability has been put forward by the authors \cite{BS1, ABS,
BSDDG}.
\begin{definition}\label{dfn:Q-net} {\bf (Q-net)}
A map $f:\bbZ^m\to\bbR\bbP^N$ is called an $m$-dimensional Q-net
(quadrilateral net, or discrete conjugate net) in $\,\bbR\bbP^N$
$(N\ge 3)$, if all its elementary quadrilaterals
$(f,f_i,f_{ij},f_j)$ (at any $u\in\bbZ^m$ and for all pairs $1\leq
i\neq j\leq m$) are planar.
\end{definition}
Thus, for any elementary quadrilateral, any representatives
$\tilde{f}$, $\tilde{f}_i$, $\tilde{f}_j$, $\tilde{f}_{ij}$ of its
vertices in the space $\bbR^{N+1}$ of homogeneous coordinates
satisfy an equation of the type
\begin{equation}\label{eq: dcn hom}
\tilde{f}_{ij}=c_{ij}\tilde{f}_j+c_{ji}\tilde{f}_i+\rho_{ij}\tilde{f}.
\end{equation}
Representatives in any hyperplane of $\bbR^{N+1}$, for instance,
in the affine part $\bbR^N$ of $\bbR\bbP^N=\bbP(\bbR^{N+1})$,
satisfy such an equation with $1=c_{ij}+c_{ji}+\rho_{ij}$, that
is,
\begin{equation} \label{eq: dcn property}
\tilde{f}_{ij}-\tilde{f}=c_{ij}(\tilde{f}_j-\tilde{f})+
c_{ji}(\tilde{f}_i-\tilde{f}).
\end{equation}

Given three points $f$, $f_1$, $f_2$ in $\bbR\bbP^N$, one can take
any point of the plane through these three points as the fourth
vertex $f_{12}$ of an elementary quadrilateral
$(f,f_1,f_{12},f_2)$ of a Q-net. Correspondingly, given any two
discrete curves $f:\bbZ\times\{0\}\to\bbR\bbP^N$ and
$f:\{0\}\times\bbZ\to\bbR\bbP^N$ with a common point $f(0,0)$, one
can construct infinitely many Q-surfaces $f:\bbZ^2\to\bbR\bbP^N$
with these curves as coordinate ones: the construction goes
inductively, on each step one has a freedom of choosing a point in
a plane (two real parameters).

On the other hand, constructing elementary hexahedra of Q-nets
corresponding to elementary 3D cubes of the lattice $\bbZ^m$
admits a well-posed initial value problem with a unique solution,
therefore one says that Q-nets are described by a {\em discrete 3D
system}:
\begin{theorem}\label{Thm: Q-cube}
{\bf (Elementary hexahedron of a Q-net)} Given seven points $f$,
$f_i$ and $f_{ij}$ $\,(1\leq i<j\leq 3)$ in $\bbR\bbP^N$, such
that each of the three quadrilaterals $(f,f_i,f_{ij},f_j)$ is
planar (i.e., $f_{ij}$ lies in the plane $\Pi_{ij}$ through $f$,
$f_i$, $f_j$), define three planes $\tau_k\Pi_{ij}$ as those
passing through the point triples $f_k$, $f_{ik}$, $f_{jk}$,
respectively. Then these three planes intersect generically at one
point:
\[
f_{123}=\tau_1\Pi_{23}\cap\tau_2\Pi_{13}\cap\tau_3\Pi_{12}\,.
\]
\end{theorem}
{\bf Proof.} Planarity of the quadrilaterals $(f,f_i,f_{ij},f_j)$
assures that all seven initial points $f$, $f_i$ and $f_{ij}$
belong to the three-dimensional space $\Pi_{123}$ through the four
points $f$, $f_1$, $f_2$, $f_3$. Hence, the planes
$\tau_k\Pi_{ij}$ lie in this three-dimensional space, and
therefore generically they intersect at exactly one point. $\Box$

\begin{figure}[htbp]
\begin{center}
\setlength{\unitlength}{0.045em}
\begin{picture}(200,240)(0,0)
 \put(0,0){\circle*{15}}    \put(150,0){\circle*{15}}
 \put(0,150){\circle*{15}}  \put(150,150){\circle*{15}}
 \put(50,50){\circle*{15}} \put(50,200){\circle*{15}}
 \put(200,50){\circle*{15}}
 \put(200,200){\circle{15}}
 \path(0,0)(150,0)       \path(0,0)(0,150)
 \path(150,0)(150,150)   \path(0,150)(150,150)
 \path(0,150)(50,200)    \path(150,150)(194,194)
 \path(50,200)(192.5,200)
 \path(200,192.5)(200,50) \path(200,50)(150,0)
 \dashline[+30]{10}(0,0)(50,50)
 \dashline[+30]{10}(50,50)(50,200)
 \dashline[+30]{10}(50,50)(200,50)
 \put(-30,-5){$f$}
 \put(-35,145){$f_3$} \put(215,45){$f_{12}$}
 \put(165,-5){$f_1$} \put(160,140){$f_{13}$}
 \put(15,50){$f_2$}  \put(10,205){$f_{23}$}
 \put(215,200){$f_{123}$}
\end{picture}
\caption{3D system on an elementary cube}\label{Fig:cube eq}
\end{center}
\end{figure}
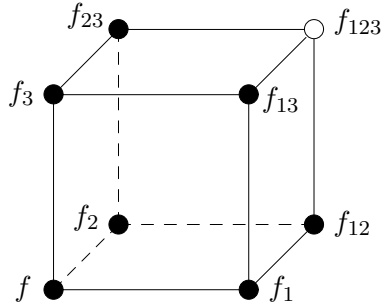

The elementary construction step from Theorem \ref{Thm: Q-cube} is
symbolically represented on Fig.~\ref{Fig:cube eq}, which is the
picture we have in mind when thinking and speaking about discrete
three-dimensional systems with dependent variables (fields)
attached to the vertices of a regular cubic lattice.
\smallskip

As follows from Theorem \ref{Thm: Q-cube}, a three-dimensional
Q-net $f:\bbZ^3\to\bbR\bbP^N$ is completely determined by its
three coordinate surfaces
\[
f:\bbZ^2\times\{0\}\to\bbR\bbP^N, \quad f:\bbZ\times\{0\}\times
\bbZ\to\bbR\bbP^N, \quad  f:\{0\}\times\bbZ^2\to\bbR\bbP^N.
\]

\begin{figure}[htbp]
\begin{center}
\setlength{\unitlength}{0.07em}
\begin{picture}(200,220)(-100,-90)

 \drawline(15,-20)(50,0)(50,47)
 \drawline(47,50)(0,50)(-35,30)(-35,-20)(15,-20)(15,30)(47,48.5)
 \drawline(15,30)(-35,30)
 \dashline{4}(-35,-20)(0,0)(0,50)\dashline{4}(0,0)(50,0)
 \drawline(30,-90)(131,-32)
 \drawline(135,-26)(135,116)
 \drawline(131,120)(-11,120)
 \drawline(-19,118)(-120,60)(-120,-90)(30,-90)(30,56)
 \drawline(34,62)(131,118)
 \drawline(26,60)(-120,60)
 \dashline{4}(-120,-90)(-15,-30)(-15,116)
 \dashline{4}(-15,-30)(131,-30)
  \dashline{2}(0,0)(-15,-30)
  \dashline{2}(-35,-20)(-120,-90)
  \dashline{2}(50,0)(132,-29)
  \dashline{2}(0,50)(-14,116)
  \dashline{2}(15,-20)(30,-90)
  \dashline{2}(-35,30)(-120,60)
  \dashline{2}(53,51.5)(131,117.5)
  \dashline{2}(15,30)(26,56)

  \put(-35,-20){\circle*{8}}       
  \put(15,-20){\circle*{8}}        
  \put(0,0){\circle*{8}}           
  \put(-35,30){\circle*{8}}        
  \put(50,0){\circle*{8}}          
  \put(0,50){\circle*{8}}          
  \put(15,30){\circle*{8}}         
  \put(50,50){\circle{8}}          

  \put(-120,-90){\circle*{10}}     
  \put(-15,-30){\circle*{9}}       
  \put(30,-90){\circle*{10}}       
  \put(-120,60){\circle*{10}}      
  \put(135,-30){\circle{10}}       
  \put(-15,120){\circle{10}}       
  \put(30,60){\circle{10}}         
  \put(130,115){$\Box$}            

  \put(-45,-15){$f$}
  \put(0,-15){$f_1$}
  \put(-18,6){$f_2$}
  \put(-50,20){$f_3$}
  \put(58,3){$f_{12}$}
  \put(-5,20){$f_{13}$}
  \put(-25,50){$f_{23}$}
  \put(58,44){$f_{123}$}
  \put(-143,-84){$f_4$}
  \put(5,-83){$f_{14}$}
  \put(-12,-41){$f_{24}$}
  \put(-143,49){$f_{34}$}
  \put(145,-28){$f_{124}$}
  \put(10,71){$f_{134}$}
  \put(-9,108){$f_{234}$}
  \put(140,103){$f_{1234}$}
\end{picture}
\caption{4D consistency of 3D systems}\label{Fig: hypercube}
\end{center}
\end{figure}
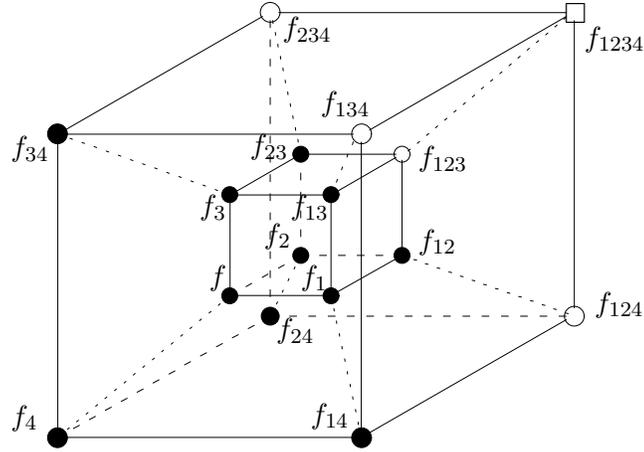

Turning to an elementary cube of the dimension $m\ge 4$, we see
that one can prescribe all points $f$, $f_i$ and $f_{ij}$ for all
$1\leq i<j\leq m$. Indeed, these data are clearly independent, and
one can construct all other vertices of an elementary cube
starting from these data, {\em provided one does not encounter
contradictions}. To see the possible source of contradictions,
consider in detail first the case of $m=4$. From $f$, $f_i$ and
$f_{ij}$ $\,(1\le i<j\le 4)$ one determines all $f_{ijk}$
uniquely. After that, one has, in principle, four different ways
to determine $f_{1234}$, from four 3D cubic faces adjacent to this
point; see Fig.~\ref{Fig: hypercube}. Absence of contradictions
means that these four values for $f_{1234}$ automatically
coincide. We call this property the 4D consistency.

\begin{definition}
{\bf (4D consistency)} A 3D system is called 4D consistent, if it
can be imposed on all three-dimensional faces of an elementary
cube of $\,\bbZ^4$.
\end{definition}

Remarkably, construction of Q-nets based on the planarity of all
elementary quadrilaterals enjoys this property.
\begin{theorem}\label{prop: dcn consistency}
{\bf (Q-nets are 4D consistent)} The 3D system governing Q-nets is
4D-consistent.
\end{theorem}
{\bf Proof.} In the construction above, the four values in
question are
\[
f_{1234}=\tau_1\tau_2\Pi_{34}\cap\tau_1\tau_3\Pi_{24}\cap
\tau_1\tau_4\Pi_{23}\,,
\]
and three other ones obtained by cyclic shifts of indices. Thus,
we have to prove that the six planes $\tau_i\tau_j\Pi_{k\ell}$
intersect in one point.

First, assume that the ambient space $\bbR\bbP^N$ has dimension
$N\ge 4$. Then, in general position, the space $\Pi_{1234}$
through the five points $f$, $f_i$ $(1\leq i\leq 4)$ is
four-dimensional. It is easy to understand that the plane
$\tau_i\tau_j\Pi_{k\ell}$ is the intersection of two
three-dimensional subspaces $\tau_i\Pi_{jk\ell}$ and
$\tau_j\Pi_{ik\ell}$. Indeed, the subspace $\tau_i\Pi_{jk\ell}$
through the four points $f_i$, $f_{ij}$, $f_{ik}$, $f_{i\ell}$
contains also $f_{ijk}$, $f_{ij\ell}$, and $f_{jk\ell}$.
Therefore, both $\tau_i\Pi_{jk\ell}$ and $\tau_j\Pi_{ik\ell}$
contain three points $f_{ij}$, $f_{ijk}$ and $f_{ij\ell}$, which
determine the plane $\tau_i\tau_j\Pi_{k\ell}$. Now the
intersection in question can be alternatively described as the
intersection of the four three-dimensional subspaces
$\tau_1\Pi_{234}$, $\tau_2\Pi_{134}$, $\tau_3\Pi_{124}$ and
$\tau_4\Pi_{123}$ of one and the same four-dimensional space
$\Pi_{1234}$. This intersection consists in the generic case of
exactly one point.

In the case of $N=3$, we embed the ambient space into
$\bbR\bbP^4$, then slightly perturb the point $f_4$ by adding a
small component in the fourth coordinate direction, then apply the
above argument, and after that send the perturbation to zero. This
proof works since, as one can easily see, on each step of the
construction the perturbation remains regular. $\Box$
\medskip

The $m$-dimensional consistency of a 3D system for $m>4$ is
defined analogously to the $m=4$ case. Remarkably and quite
generally, the 4-dimensional consistency already implies
$m$-dimensional consistency for all $m>4$.
\begin{theorem}\label{Thm: 4D yields MD}
{\bf (4D consistency yields consistency in all higher dimensions)}
Any 4D consistent discrete 3D system is also $m$-dimensionally
consistent for any $m>4$.
\end{theorem}
{\bf Proof} goes by induction from the $(m-1)$-dimensional
consistency to the $m$-dimensional consistency, but, for the sake
of notational simplicity, we present the details for the case
$m=5$ only, the general case being absolutely similar.

Initial data for a 3D system on the 5D cube $\cC_{12345}$ with the
fields on vertices consist of the fields $f$, $f_i$ and $f_{ij}$
for all $1\le i<j\le 5$. From these data one first gets ten fields
$f_{ijk}$ for $1\le i<j<k\le 5$, and then five fields
$f_{ijk\ell}$ for $1\le i<j<k<\ell\le 5$ (the fact that the latter
are well defined is nothing but the assumed 4D consistency for the
4D cubes $\cC_{ijkl}$). Now, one has ten possibly different values
for $f_{12345}$, coming from ten 3D cubes $\tau_i\tau_j\cC_{k\ell
m}$. To prove that these ten values coincide, consider five 4D
cubes $\tau_i\cC_{jk\ell m}$. For instance, for the 4D cube
$\tau_1\cC_{2345}$ the assumed consistency assures that the four
values for $f_{12345}$ coming from four 3D cubes
\[
\tau_1\tau_2\cC_{345},\quad \tau_1\tau_3\cC_{245},\quad
\tau_1\tau_4\cC_{235},\quad \tau_1\tau_5\cC_{234}
\]
are all the same. Similarly, for the 4D cube $\tau_2\cC_{1345}$
the 4D consistency leads to the conclusion that the four values
for $f_{12345}$ coming from
\[
\tau_1\tau_2\cC_{345},\quad \tau_2\tau_3\cC_{145},\quad
\tau_2\tau_4\cC_{135},\quad \tau_2\tau_5\cC_{134}
\]
coincide. Note that the 3D cube $\tau_1\tau_2\cC_{345}$, the
intersection of $\tau_1\cC_{2345}$ and $\tau_2\cC_{1345}$, is
present in both lists, so that we now have seven coinciding values
for $f_{12345}$. Adding similar conclusions for other 4D cubes
$\tau_i\cC_{jk\ell m}$, we arrive at the desired result. $\Box$
\smallskip

Theorems \ref{prop: dcn consistency}, \ref{Thm: 4D yields MD}
yield that Q-nets are $m$-dimensionally consistent for any $m\ge
4$. This fact, in turn, yields the existence of transformations of
Q-nets with remarkable permutability properties. Referring for
details to \cite{DSM, BSDDG}, we mention here only the definition.
\begin{definition}\label{Def: Jonas}
{\bf (F-transformation of Q-nets)} Two $m$-dimensional Q-nets
$f,f^+:\bbZ^m\to\bbR\bbP^N$ are called F-transforms (fundamental
transforms) of one another, if all quadrilaterals
$(f,f_i,f_i^+,f^+)$ (at any $u\in\bbZ^m$ and for all $\,1\le i\le
m$) are planar, i.e., if the net
$F:\bbZ^m\times\{0,1\}\to\bbR\bbP^N$ defined by $F(u,0)=f(u)$ and
$F(u,1)=f^+(u)$ is a two-layer $(m+1)$-dimensional Q-net.
\end{definition}
It follows from Theorem \ref{Thm: Q-cube} that, given a Q-net $f$,
its F-transform $f^+$ is uniquely defined as soon as its points
along the coordinate axes are suitably prescribed.

\subsection{Discrete line congruences}
\label{subsect: congr}

Another important geometrical objects described by a discrete 3D
system which is 4D consistent, are {\em discrete line
congruences}. Their theory has been developed by Doliwa, Santini
and Ma\~{n}as \cite{DSM}, whose presentation we follow in this
section.

Let $\cL^N$ be the space of lines in $\bbR\bbP^N$; it can be
identified with the Grassmannian ${\rm Gr}(N+1,2)$ of
two-dimensional vector subspaces of $\bbR^{N+1}$.

\begin{definition}\label{Def: congr}
{\bf (Discrete line congruence)} A map $\ell:\bbZ^m\to\cL^N$ is
called an $m$-dimensional discrete line congruence in $\bbR\bbP^N$
$(N\ge 3)$, if any two neighboring lines $\ell$, $\ell_i$ (at any
$u\in\bbZ^m$ and for any $1\le i\le m$) intersect (are co-planar).
\end{definition}

For instance, lines $\ell=(ff^+)$ connecting corresponding points
of two Q-nets $f,f^+:\bbZ^m\to\bbR\bbP^N$ in the relation of
F-transformation clearly build a discrete line congruence.

A discrete line congruence is called {\em generic}, if for any
$u\in\bbZ^m$ and for any $1\le i\neq j\neq k\neq i\le m$, the four
lines $\ell$, $\ell_i$, $\ell_j$ and $\ell_k$ span a
four-dimensional space (i.e., a space of a maximal possible
dimension). This yields, in particular, that for any $u\in\bbZ^m$
and for any $1\le i\neq j\le m$, the three lines $\ell$, $\ell_i$
and $\ell_j$ span a three-dimensional space.

Construction of line congruences is similar to that of Q-nets.
Given three lines $\ell$, $\ell_1$, $\ell_2$ of a congruence, one
has a two-parameter family of lines admissible as the fourth one
$\ell_{12}$: connect by a line any point of $\ell_1$ with any
point of $\ell_2$. Thus, given any two sequences of lines
$\ell:\bbZ\times\{0\}\to\cL^N$ and $\ell:\{0\}\times\bbZ\to\cL^N$
such that any two neighboring lines are co-planar, one can extend
them to a two-dimensional line congruence $f:\bbZ^2\to\cL^N$ in an
infinite number of ways: on each step of the inductive procedure
one has a freedom of choosing a line from a two-parameter family.

The next theorem shows that non-degenerate line congruences are
described by a {\em discrete 3D system}:
\begin{theorem}\label{Thm:  congr cube}
{\bf (Elementary hexahedron of a discrete line congruence)} Given
seven lines $\ell$, $\ell_i$ and $\ell_{ij}$ $\,(1\leq i<j\leq 3)$
in $\bbR\bbP^N$, such that $\ell$ intersects each of $\ell_i$, the
space $V_{123}$ spanned by $\ell$, $\ell_1$, $\ell_2$, $\ell_3$
has dimension four, and each $\ell_i$ intersects both $\ell_{ij}$
and $\ell_{ik}$, there is a unique line $\ell_{123}$ that
intersects all three $\ell_{ij}$.
\end{theorem}
{\bf Proof.} All seven lines, and therefore also the
three-dimensional spaces $\tau_iV_{jk}={\rm
span}(\ell_i,\ell_{ij},\ell_{ik})$ lie in $V_{123}$. A line that
intersects all three of $\ell_{ij}$ should lie in the intersection
of these three three-dimensional spaces. But a generic
intersection of three three-dimensional spaces in $V_{123}$ is a
line:
\[
\ell_{123}=\tau_1V_{23}\cap \tau_2V_{13}\cap\tau_3V_{12}.
\]
It is now not difficult to realize that this line does, indeed,
intersect all three of $\ell_{ij}$. For instance,
$\tau_1V_{23}\cap \tau_2V_{13}={\rm
span}(\ell_{12},\ell_{13})\cap{\rm span}(\ell_{12},\ell_{23})$ is
a plane containing $\ell_{12}$, therefore its intersection with
$\tau_3V_{12}$ (the line $\ell_{123}$) intersects $\ell_{12}$.
$\Box$
\smallskip

A similar argument shows:
\begin{theorem}\label{prop: congr consistency}
{\bf (Discrete line congruences are 4D consistent)} The 3D system
governing discrete line congruences is 4D-consistent.
\end{theorem}

Like in the case of Q-nets, this theorem yields the existence of
transformations of discrete line congruences with remarkable
permutability properties.
\begin{definition}\label{Def: XXX}
{\bf (F-transformation of line congruences)} Two $m$-dimen\-sional
line congruences $\ell,\ell^+:\bbZ^m\to\cL^N$ are called
F-transforms of one another, if the corresponding lines $\ell$ and
$\ell^+$ intersect (at any $u\in\bbZ^m$), i.e., if the map
$L:\bbZ^m\times\{0,1\}\to\cL^N$ defined by $L(u,0)=\ell(u)$ and
$L(u,1)=\ell^+(u)$ is a two-layer $(m+1)$-dimensional line
congruence.
\end{definition}
Again, it follows from Theorem \ref{Thm: Q-cube} that, given a
line congruence $\ell$, its F-transform $\ell^+$ is uniquely
defined as soon as its lines along the coordinate axes are
suitably prescribed.
\smallskip

\begin{figure}[htbp]
 \psfrag{p^i}[Bl][bl][0.9]{$f^{(i)}$}
 \psfrag{p^j}[Bl][bl][0.9]{$f^{(j)}$}
 \psfrag{p_j^i}[Bl][bl][0.9]{$f_j^{(i)}$}
 \psfrag{p_i^j}[Bl][bl][0.9]{$f_i^{(j)}$}
 \psfrag{l}[Bl][bl][0.9]{$\ell$}
 \psfrag{l_i}[Bl][bl][0.9]{$\ell_i$}
 \psfrag{l_j}[Bl][bl][0.9]{$\ell_j$}
 \psfrag{l_ij}[Bl][bl][0.9]{$\ell_{ij}$}
 \center{\includegraphics[width=120mm]{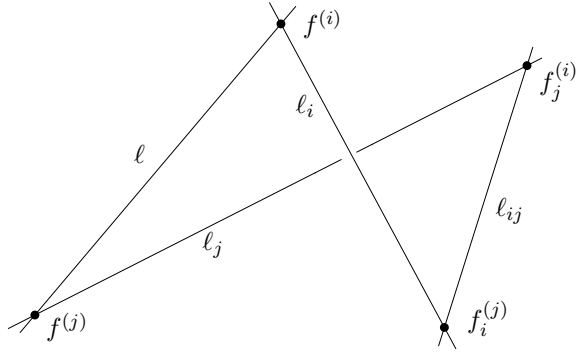}}
 \caption{Four lines of a congruence}
 \label{Fig: congr}
\end{figure}

According to Definition \ref{Def: congr}, any two neighboring
lines $\ell=\ell(u)$ and $\ell_i=\ell(u+e_i)$ of a line congruence
intersect at exactly one point $f=\ell\cap\ell_i\in\bbR\bbP^N$
which is thus combinatorially associated with the edge $(u,u+e_i)$
of the lattice $\bbZ^m$: $f=f(u,u+e_i)$. It is, however, sometimes
more convenient to use the notation $f(u,u+e_i)=f^{(i)}(u)$ for
this points, thus associating it to the vertex $u$ of the lattice
(and, of course, to the coordinate direction $i$). See
Fig.~\ref{Fig: congr}.
\begin{definition}\label{Def: focal net}
{\bf (Focal net)} For a discrete line congruence
$\ell:\bbZ^m\to\cL^N$, the map $f^{(i)}:\bbZ^m\to\bbR\bbP^N$
defined by $f^{(i)}(u)=\ell(u)\cap\ell(u+e_i)$ is called its
$i$-th focal net.
\end{definition}

\begin{theorem}\label{thm: focal Q-nets}
For a non-degenerate discrete line congruence
$\ell:\bbZ^m\to\cL^N$, all its focal nets
$f^{(k)}:\bbZ^m\to\bbR\bbP^N$, $1\le k\le m$, are Q-nets.
\end{theorem}
{\bf Proof} consists of two steps.

$\blacktriangleright\;$ First, one shows that for the $k$-th focal
net $f^{(k)}$, all elementary quadrilaterals
$(f^{(k)},f_i^{(k)},f_{ik}^{(k)},f_k^{(k)})$ are planar. This is
true for any line congruence. Indeed, both points $f^{(k)}$ and
$f_k^{(k)}$ lie on the line $\ell_k$, while both points
$f_i^{(k)}$ and $f_{ik}^{(k)}$ lie on the line $\ell_{ik}$.
Therefore, all four points lie in the plane spanned by these two
lines $\ell_k$ and $\ell_{ik}$ which intersect by definition of a
line congruence.

\begin{figure}[htbp]
 \psfrag{p^k}[Bl][bl][0.9]{$f^{(k)}$}
 \psfrag{p_i^k}[Bl][bl][0.9]{$f_i^{(k)}$}
 \psfrag{p_j^k}[Bl][bl][0.9]{$f_j^{(k)}$}
 \psfrag{p_ij^k}[Bl][bl][0.9]{$f_{ij}^{(k)}$}
 \psfrag{l}[Bl][bl][0.9]{$\ell$}
 \psfrag{l_i}[Bl][bl][0.9]{$\ell_i$}
 \psfrag{l_j}[Bl][bl][0.9]{$\ell_j$}
 \psfrag{l_k}[Bl][bl][0.9]{$\ell_k$}
 \psfrag{l_ij}[Bl][bl][0.9]{$\ell_{ij}$}
 \psfrag{l_ik}[Bl][bl][0.9]{$\ell_{ik}$}
 \psfrag{l_jk}[Bl][bl][0.9]{$\ell_{jk}$}
 \psfrag{l_ijk}[Bl][bl][0.9]{$\ell_{ijk}$}
 \center{\includegraphics[width=140mm]{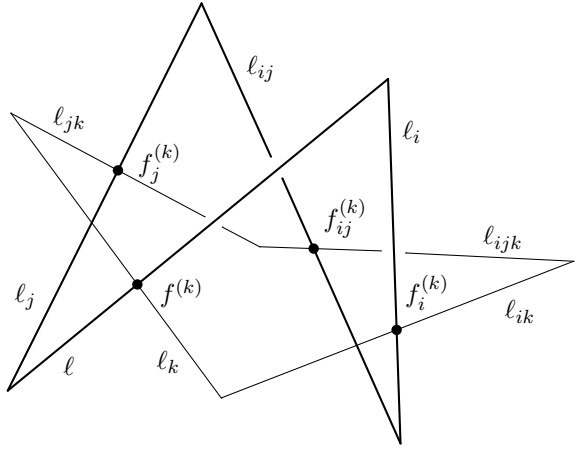}}
 \caption{Elementary $(ij)$ quadrilateral of the $k$-th focal net}
 \label{Fig: Focal quad}
\end{figure}

$\blacktriangleright\;$ Second, one shows that for the $k$-th
focal net $f^{(k)}$, all elementary quadrilaterals
$(f^{(k)},f_i^{(k)},f_{ij}^{(k)},f_j^{(k)})$, with both $i\neq j$
different from $k$, are planar. Here, one uses essentially the
assumption that the line congruence $\ell$ is generic. All four
points in question lie in each of the three-dimensional spaces
\[
V_{ij}={\rm span}(\ell,\ell_i,\ell_j,\ell_{ij}) \quad{\rm
and}\quad \tau_kV_{ij}={\rm
span}(\ell_k,\ell_{ik},\ell_{jk},\ell_{ijk})
\]
(see Fig.~\ref{Fig: Focal quad}). Both 3-spaces lie in the
four-dimensional space $V_{ijk}={\rm
span}(\ell,\ell_i,\ell_j,\ell_k)$, so that generically their
intersection is a plane. $\Box$

\begin{corollary}
{\bf (Focal net of F-transformation of a line congruence)}  Given
two generic line congruences $\ell,\ell^+:\bbZ^m\to\cL^N$ in the
relation of F-transformation, the intersection points
$f=\ell\cap\ell^+$ form a Q-net $f:\bbZ^m\to\bbR\bbP^N$.
\end{corollary}

\subsection{Q-nets in quadrics}
 \label{Sect: dcn quadr}

We consider an important admissible reduction of Q-nets: they can
be consistently restricted to an arbitrary quadric in
$\bbR\bbP^N$. In the smooth differential geometry, i.e., for
conjugate nets, this is due to Darboux \cite{Da1}. In the discrete
differential geometry this result has been found by Doliwa
\cite{D2}.

A deep reason for this result is the following fundamental fact
well known in classical projective geometry (see, e.g.,
\cite{Bl}):
\begin{theorem}\label{Th: associated point}
{\bf (Associated point)}  For any seven points of $\,\bbC\bbP^3$
in general position, there exists the eighth point (called the
associated one), which belongs to any quadric through the original
seven points.
\end{theorem}
{\bf Proof} is based on the following computations. The equation
$\cQ=0$ of a quadric in $\bbC\bbP^3$ has ten coefficients
(homogeneous polynomial of 4 variables). Therefore, a unique
quadric $\cQ=0$ can be drawn through nine points in general
position. Similarly, a pencil (one-parameter linear family) of
quadrics $\cQ+\lambda\cQ'=0$ can be drawn through eight points in
general position, and a two-parameter linear family of quadrics
$\cQ+\lambda\cQ'+\mu\cQ''=0$ can be drawn through seven points in
general position. Generically, solution of a system of three
quadratic equations
\[
\cQ=0,\quad \cQ'=0,\quad \cQ''=0
\]
for the intersection of three quadrics in $\bbC\bbP^3$ consists of
eight points. It can be shown that the three quadrics spanning the
above-mentioned two-parameter family can be considered generic
enough for such a conclusion. Clearly, the resulting eight points
lie on every quadric of the two-parameter family. $\Box$

\begin{theorem}\label{Th: Q-net in Q}
{\bf (Elementary hexahedron of a Q-net in a quadric)} If seven
points $f$, $f_i$, and $f_{ij}\,$ ($1\le i<j\le 3$) of an
elementary hexahedron of a Q-net $f:\bbZ^m\to\bbR\bbP^N$ belong to
a quadric $\,\cQ\subset\bbR\bbP^N$, then so does the eighth point
$f_{123}$.
\end{theorem}
{\bf Proof.} The original seven points can be considered lying in
a three-dimen\-sional space, and they are known to belong to three
(degenerate) quadrics -- the pairs of planes
$\Pi_{jk}\cup\tau_i\Pi_{jk}$ for $(jk)=(12), (23), (31)$. Clearly,
the eighth intersection point of these quadrics is
$f_{123}=\tau_1\Pi_{23}\cap\tau_2\Pi_{31}\cap\tau_3\Pi_{12}$, and
this has to be the associated point. According to Theorem \ref{Th:
associated point}, it belongs to any quadric through the original
seven points, in particular, to $\cQ$. $\Box$

\section{Geometries of spheres}
\label{Sect: geom}

\subsection{Lie geometry}
\label{subsect: Lie}

A classical source on {\em Lie geometry} is Blaschke's book
\cite{BlIII}, see also a modern account by Cecil \cite{Ce}.

\begin{figure}[htbp]
 \psfrag{x}[Bl][bl][0.9]{$x$}
 \psfrag{P}[Bl][bl][0.9]{$P$}
 \psfrag{v}[Bl][bl][0.9]{$v$}
 \center{\includegraphics[width=120mm]{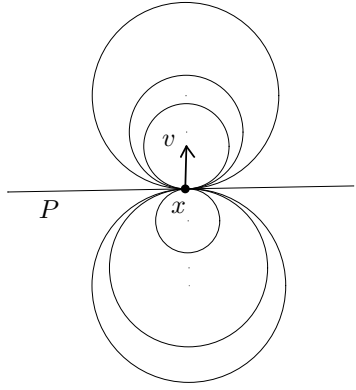}}
 \caption{Contact element}
 \label{Fig: contact}
\end{figure}

Following geometric objects in the Euclidean space $\bbR^N$ are
elements of Lie geometry:
\begin{itemize}
\item {\em Oriented hyperspheres.} A hypersphere in $\bbR^N$ with
center $c\in\bbR^N$ and radius $r>0$ is described by the equation
$S=\{x\in\bbR^N: |x-c|^2=r^2\}$. It divides $\bbR^N$ in two parts,
inner and outer. If one denotes one of two parts of $\bbR^N$ as
positive, one comes to the notion of an oriented hypersphere.
Thus, there are two oriented hyperspheres $S^{\pm}$ for any $S$.
One can take the orientation of a hypersphere into account by
assigning a {\em signed radius} $\pm r$ to it. For instance, one
can assign positive radii $r>0$ to hyperspheres with the inward
field of unit normals and negative radii $r<0$ to hyperspheres
with the outward field of unit normals.

\item {\em Oriented hyperplanes.} A hyperplane in $\bbR^N$ is
given by the equation $P=\{x\in\bbR^N: \langle v,x\rangle=d\}$,
with a unit normal $v\in\bbS^{N-1}$ and $d\in\bbR$. Clearly, the
pairs $(v,d)$ and $(-v,-d)$ represent one and the same hyperplane.
It divides $\bbR^N$ in two halfspaces. Denoting one of two
halfspaces as positive, one arrives at the notion of an oriented
hyperplane. Thus, there are two oriented hyperplanes $P^{\pm}$ for
any $P$. One can take the orientation of a hypersphere into
account by assigning the pair $(v,d)$ to the hyperplane with the
unit normal $v$ pointing into the positive halfspace.

\item{\em Points.} One considers points $x\in\bbR^N$ as
hyperspheres of a vanishing radius.

\item{\em Infinity.} One compactifies the space $\bbR^N$ by adding
the point at infinity $\infty$, with the understanding that a
basis of open neighborhoods of $\infty$ is given, e.g., by the
outer parts of the hyperspheres $|x|^2=r^2$. Topologically the so
defined compactification is equivalent to a sphere $\bbS^N$.

\item {\em Contact elements.} A contact element of a hypersurface
is a pair consisting of a point $x\in\bbR^N$ and an (oriented)
hyperplane $P$ through $x$; alternatively, one can use a normal
vector $v$ to $P$ at $x$. In the framework of Lie geometry, a
contact element can be identified with a set (a pencil) of all
hyperspheres $S$ through $x$ which are in an oriented contact with
$P$ (and with one another), thus sharing the normal vector $v$ at
$x$, see Fig.~\ref{Fig: contact}.
\end{itemize}

All these elements are modelled in Lie geometry as points, resp.
lines, in the $(N+2)$-dimensional projective space
$\bbP(\bbR^{N+1,2})$ with the space of homogeneous coordinates
$\bbR^{N+1,2}$. The latter is the space spanned by $N+3$ linearly
independent vectors $\ee_1,\ldots,\ee_{N+3}$ and equipped with the
pseudo-euclidean scalar product
\[
\langle \ee_i,\ee_j\rangle=\left\{\begin{array}{rl} 1, &
i=j\in\{1,\ldots, N+1\},\\ -1, & i=j\in\{N+2,N+3\},\\ 0, & i\neq
j.\end{array}\right.
\]
It is convenient to introduce two isotropic vectors
\begin{equation}\label{eq: e0infty}
\ee_0=\tfrac{1}{2}(\ee_{N+2}-\ee_{N+1}),\quad
\ee_\infty=\tfrac{1}{2}(\ee_{N+2}+\ee_{N+1}),
\end{equation}
for which
\[
\langle \ee_0,\ee_0\rangle=\langle
\ee_\infty,\ee_\infty\rangle=0,\quad \langle
\ee_0,\ee_\infty\rangle=-\tfrac{1}{2}.
\]
The models of the above elements in the space $\bbR^{N+1,2}$ of
homogeneous coordinates are as follows:
\begin{itemize}
\item {\em Oriented hypersphere with center $c\in\bbR^N$ and
signed radius $r\in\bbR$:}
\begin{equation}\label{eq: Lie sphere}
\hat{s}=c+\ee_0+(|c|^2-r^2)\ee_\infty+r\ee_{N+3}.
\end{equation}
\item {\em Oriented hyperplane $\langle v,x\rangle=d$ with
$v\in\bbS^{N-1}$ and $d\in\bbR$}:
\begin{equation}\label{eq: Lie plane}
\hat{p}=v+0\cdot\ee_0+2d\ee_\infty+\ee_{N+3}.
\end{equation}
\item {\em Point $x\in\bbR^N$:}
\begin{equation}\label{eq: Lie point}
\hat{x}=x+\ee_0+|x|^2\ee_\infty+0\cdot\ee_{N+3}.
\end{equation}
\item {\em Infinity $\infty$:}
\begin{equation}\label{eq: Lie infty}
\hat{\infty}=\ee_\infty.
\end{equation}
\item {\em Contact element $(x,P)$:}
\begin{equation}\label{eq: Lie contact el}
{\rm span}(\hat{x},\hat{p}).
\end{equation}
\end{itemize}
In the projective space $\bbP(\bbR^{N+1,2})$ the first four types
of elements are represented by the points which are equivalence
classes of (\ref{eq: Lie sphere})--(\ref{eq: Lie infty}) with
respect to the relation $\xi\sim\eta\Leftrightarrow
\xi=\lambda\eta$ with $\lambda\in\bbR^*$ for
$\xi,\eta\in\bbR^{N+1,2}$. A contact element is represented by the
line in $\bbP(\bbR^{N+1,2})$ through the points with the
representatives $\hat{x}$ and $\hat{p}$. We mention several
fundamentally important features of this model:
\begin{itemize}
\item[(i)] All the above elements belong to the {\em Lie quadric}
$\bbP(\bbL^{N+1,2})$, where
\begin{equation}\label{eq: Lie quadric}
\bbL^{N+1,2}=\big\{\xi\in\bbR^{N+1,2}: \langle
\xi,\xi\rangle=0\big\}.
\end{equation}
Moreover, points of $\bbP(\bbL^{N+1,2})$ are in a one-to-one
correspondence with oriented hyperspheres in $\bbR^N$, including
degenerate case: proper hyperspheres correspond to points of
$\bbP(\bbL^{N+1,2})$ with both $\ee_0$- and $\ee_{N+3}$-components
non-vanishing, hyperplanes correspond to points of
$\bbP(\bbL^{N+1,2})$ with vanishing $\ee_0$-component, points
correspond to points of $\bbP(\bbL^{N+1,2})$ with vanishing
$\ee_{N+3}$-component, and infinity corresponds to the only point
of $\bbP(\bbL^{N+1,2})$ with both $\ee_0$- and
$\ee_{N+3}$-components vanishing.

\item[(ii)] Two oriented hyperspheres $S_1,S_2$ are in an oriented
contact (i.e., are tangent to each other with the unit normals at
tangency pointing in the same direction), if and only if
\begin{equation}\label{eq: sph tang}
|c_1-c_2|^2=(r_1-r_2)^2,
\end{equation}
and this is equivalent to $\langle \hat{s}_1,\hat{s}_2\rangle=0$.

\item[(iii)] An oriented hypersphere $S=\{x\in\bbR^N:
|x-c|^2=r^2\}$ is in an oriented contact with an oriented
hyperplane $P=\{x\in\bbR^N: \langle v,x\rangle=d\}$, if and only
if
\begin{equation}\label{eq: tangency}
\langle c,v\rangle-r-d=0.
\end{equation}
Indeed, equation of the hyperplane $P$ tangent to $S$ at $x_0\in
S$ reads: $\langle x_0-c,x-c\rangle=r^2$. Denoting by
$v=(c-x_0)/r$ the unit normal vector of $P$ (recall that the
positive radii are assigned to spheres with inward unit normals),
we can write the above equation as $\langle v,x\rangle=d$ with
$d=\langle c,(c-x_0)/r\rangle-r=\langle c,v\rangle-r$, which
proves eq. (\ref{eq: tangency}). Now, the latter equation is
equivalent to $\langle \hat{s},\hat{p}\rangle=0$.
\smallskip

\item[(iv)] A point $x$ can be considered as a hypersphere of
radius $r=0$ (in this case both oriented hyperspheres coincide).
An incidence relation $x\in S$ with a hypersphere $S$ (resp. $x\in
P$ with a hyperplane $P$) can be interpreted as a particular case
of an oriented contact of a sphere of radius $r=0$ with $S$ (resp.
with $P$), and it takes place if and only if
$\langle\hat{x},\hat{s}\rangle=0$ (resp.
$\langle\hat{x},\hat{p}\rangle=0$).

\item[(v)] For any hyperplane $P$, there holds
$\langle\hat{\infty},\hat{p}\rangle=0$. One can interpret
hyperplanes as hyperspheres (of an infinite radius) through
$\infty$. More precisely, a hyperplane $\langle v,x\rangle=d$ can
be interpreted as a limit, as $r\to\infty$, of the hyperspheres of
radii $r$ with the centers located at $c=rv+u$, with $\langle
v,u\rangle=d$. Indeed, the representatives (\ref{eq: Lie sphere})
of such spheres are
\begin{eqnarray*}
\hat{s} & = & (rv+u)+\ee_0+(2dr+\langle u,u\rangle)\ee_\infty+r\ee_{N+3}\\
        & \sim &
        (v+O(1/r))+(1/r)\ee_0+(2d+O(1/r))\ee_\infty+\ee_{N+3}\\
        & = & \hat{p}+O(1/r).
\end{eqnarray*}
Moreover, for similar reasons, the infinity $\infty$ can be
considered as a limiting position of any sequence of points $x$
with $|x|\to\infty$.

\item[(vi)] Any two hyperspheres $S_1$, $S_2$ in an oriented
contact determine a contact element (their point of contact and
their common tangent hyperplane). For their representatives
$\hat{s}_1$, $\hat{s}_2$ in $\bbR^{N+1,2}$, the line in
$\bbP(\bbR^{N+1,2})$ through the corresponding points in
$\bbP(\bbL^{N+1,2})$ is {\em isotropic}, i.e., lies entirely on
the Lie quadric $\bbP(\bbL^{N+1,2})$. This follows from
\[
\langle\alpha_1\hat{s}_1+\alpha_2\hat{s}_2,
\alpha_1\hat{s}_1+\alpha_2\hat{s}_2\rangle=
2\alpha_1\alpha_2\langle\hat{s}_1,\hat{s}_2\rangle=0.
\]
Such a line contains exactly one point whose representative
$\hat{x}$ has vanishing $\ee_{N+3}$-component (and corresponds to
$x$, the common point of contact of all hyperspheres), and, if
$x\neq\infty$, exactly one point whose representative $\hat{p}$
has vanishing $\ee_0$-component (and corresponds to $P$, the
common tangent hyperplane of all hyperspheres). In case when an
isotropic line contains $\hat{\infty}$, all its points represent
parallel hyperplanes, which constitute a contact element through
$\infty$.

\end{itemize}

Thus, if one considers hyperplanes as hyperspheres of infinite
radii, and points as hyperspheres of vanishing radii, then one can
conclude that:
\begin{itemize}
\item[$\blacktriangleright$] oriented hyperspheres are in a
one-to-one correspondence with points of the Lie quadric
$\bbP(\bbL^{N+1,2})$ in the projective space $\bbP(\bbR^{N+1,2})$;

\item[$\blacktriangleright$] oriented contact of two oriented
hyperspheres corresponds to orthogonality of (any) representatives
of the corresponding points in $\bbP(\bbR^{N+1,2})$.

\item[$\blacktriangleright$] contact elements of hypersurfaces are
in a one-to-one correspondence with isotropic lines in
$\bbP(\bbR^{N+1,2})$. We will denote the set of all such lines by
$\cL_0^{N+1,2}$.
\end{itemize}

According to F.~Klein's Erlangen Program, Lie geometry is the
study of properties of transformations which map oriented
hyperspheres (including points and hyperplanes) to oriented
hyperspheres and, moreover, preserve the oriented contact of
hypersphere pairs. In the projective model described above, Lie
geometry is the study of projective transformations of
$\bbP(\bbR^{N+1,2})$ which leave $\bbP(\bbL^{N+1,2})$ invariant,
and, moreover, preserve orthogonality of points of
$\bbP(\bbL^{N+1,2})$ (which is understood as orthogonality of
their lifts to $\bbL^{N+1,2}\subset\bbR^{N+1,2}$; clearly, this
relation does not depend on the choice of lifts). Such
transformations are called {\em Lie sphere transformations}.
\begin{theorem}\label{thm: Lie fund}
{\bf (Fundamental theorem of Lie geometry)}

\noindent a) The group of Lie sphere transformations is isomorphic
to $O(N+1,2)/\{\pm I\}$.

\noindent b) Every line preserving diffeomorphism of
$\,\bbP(\bbL^{N+1,2})$ is the restriction to $\bbP(\bbL^{N+1,2})$
of a Lie sphere transformation.
\end{theorem}

Since (non-)vanishing of the $\ee_0$- or of the
$\ee_{N+3}$-component of a point in $\bbP(\bbL^{N+1,2})$ is not
invariant under a general Lie sphere transformation, there is no
distinction between oriented hyperspheres, oriented hyperplanes
and points in Lie geometry.

\subsection{M\"obius geometry}
\label{subsect: Moebius}

Blaschke's book \cite{BlIII} serves also as a classical source on
{\em M\"obius geometry}, a modern account can be found in
\cite{HJ2}.

M\"obius geometry is a subgeometry of Lie geometry, with points
distinguishable among all hyperspheres as those of radius zero.
Thus, M\"obius geometry studies properties of hyperspheres
invariant under the subgroup of Lie sphere transformations
preserving the set of points. In the projective model, points of
$\bbR^N$ are distinguished as points of $\bbP(\bbL^{N+1,2})$ with
the vanishing $\ee_{N+3}$-component. (Of course, one could replace
here $\ee_{N+3}$ by any time-like vector.) Thus, M\"obius geometry
studies the subgroup of Lie sphere transformations preserving the
subset of $\bbP(\bbL^{N+1,2})$ with the vanishing
$\ee_{N+3}$-component. Following geometric objects in $\bbR^N$ are
elements of M\"obius geometry.

\begin{itemize}
\item {\em (Non-oriented) hyperspheres} $S=\{x\in\bbR^N:
|x-c|^2=r^2\}$ with centers $c\in\bbR^N$ and radii $r>0$.

\item {\em (Non-oriented) hyperplanes} $P=\{x\in\bbR^N: \langle
v,x\rangle=d\}$, with unit normals $v\in\bbS^{N-1}$ and
$d\in\bbR$.

\item{\em Points} $x\in\bbR^N$.

\item{\em Infinity} $\infty$ which compactifies $\bbR^N$ into
$\bbS^N$.
\end{itemize}

In modelling these elements, one can use the Lie-geometric
description and just omit the $\ee_{N+3}$-component. The resulting
objects are points of the $(N+1)$-dimensional projective space
$\bbP(\bbR^{N+1,1})$ with the space of homogeneous coordinates
$\bbR^{N+1,1}$. The latter is the space spanned by $N+2$ linearly
independent vectors $\ee_1,\ldots,\ee_{N+2}$ and equipped with the
Minkowski scalar product
\[
\langle \ee_i,\ee_j\rangle=\left\{\begin{array}{rl} 1, &
i=j\in\{1,\ldots, N+1\},\\ -1, & i=j=N+2,\\ 0, & i\neq
j.\end{array}\right.
\]
We continue to use notations (\ref{eq: e0infty}) in the context of
the M\"obius geometry. The above elements are modelled in the
space $\bbR^{N+1,1}$ of homogeneous coordinates as follows:
\begin{itemize}
\item {\em Hypersphere with center $c\in\bbR^N$ and radius $r>0$:}
\begin{equation}\label{eq: Moeb sphere}
\hat{s}=c+\ee_0+(|c|^2-r^2)\ee_\infty.
\end{equation}
\item {\em Hyperplane $\langle v,x\rangle=d$ with $v\in\bbS^{N-1}$
and $d\in\bbR$}:
\begin{equation}\label{eq: Moeb plane}
\hat{p}=v+0\cdot\ee_0+2d\ee_\infty.
\end{equation}
\item {\em Point $x\in\bbR^N$:}
\begin{equation}\label{eq: Moeb point}
\hat{x}=x+\ee_0+|x|^2\ee_\infty.
\end{equation}
\item {\em Infinity $\infty$:}
\begin{equation}\label{eq: Moeb infty}
\hat{\infty}=\ee_\infty.
\end{equation}
\end{itemize}
In the projective space $\bbP(\bbR^{N+1,1})$ these elements are
represented by points which are equivalence classes of (\ref{eq:
Moeb sphere})--(\ref{eq: Moeb infty}) with respect to the usual
relation $\xi\sim\eta\Leftrightarrow \xi=\lambda\eta$ with
$\lambda\in\bbR^*$ for $\xi,\eta\in\bbR^{N+1,1}$. Fundamental
features of these identifications:
\begin{itemize}
\item[(i)] The infinity $\hat{\infty}$ can be considered as a
limit of any sequence of $\hat{x}$ for $x\in\bbR^N$ with
$|x|\to\infty$. Elements $x\in\bbR^N\cup\{\infty\}$ are in a
one-to-one correspondence with points of the projectivized {\em
light cone} $\bbP(\bbL^{N+1,1})$, where
\begin{equation}\label{eq: Moeb quadric}
\bbL^{N+1,1}=\big\{\xi\in\bbR^{N+1,1}: \langle
\xi,\xi\rangle=0\big\}.
\end{equation}
Points $x\in\bbR^N$ correspond to points of $\bbP(\bbL^{N+1,1})$
with a non-vanishing $\ee_0$-component, while $\infty$ corresponds
to the only point of $\bbP(\bbL^{N+1,1})$ with the vanishing
$\ee_0$-component.

\item[(ii)] Hyperspheres $\hat{s}$ and hyperplanes $\hat{p}$
belong to $\bbP(\bbR^{N+1,1}_{\rm out})$, where
\begin{equation}\label{eq: Moeb Rout}
\bbR^{N+1,1}_{\rm out}=\big\{\xi\in\bbR^{N+1,1}: \langle
\xi,\xi\rangle>0\big\}
\end{equation}
is the set of space-like vectors of the Minkowski space
$\bbR^{N+1,1}$. Hyperplanes can be interpreted as hyperspheres (of
an infinite radius) through $\infty$.

\item[(iii)] Two hyperspheres $S_1,S_2$ with centers $c_1,c_2$ and
radii $r_1,r_2$ intersect orthogonally, if and only if
\begin{equation}\label{eq: ortho sph}
|c_1-c_2|^2=r_1^2+r_2^2,
\end{equation}
which is equivalent to $\langle\hat{s}_1,\hat{s}_2\rangle=0$.
Similarly, a hypersphere $S$ intersects orthogonally with a
hyperplane $P$, if and only if its center lies in $P$:
\begin{equation}
\langle c,v\rangle-d=0,
\end{equation}
which is equivalent to $\langle \hat{s},\hat{p}\rangle=0$.

\item[(iv)] A point $x$ can be considered as a limiting case of a
hypersphere with radius $r=0$. An incidence relation $x\in S$ with
a hypersphere $S$ (resp. $x\in P$ with a hyperplane $P$) can be
interpreted as a particular case of an orthogonal intersection of
a sphere of radius $r=0$ with $S$ (resp. with $P$), and it takes
place if and only if $\langle\hat{x},\hat{s}\rangle=0$ (resp.
$\langle\hat{x},\hat{p}\rangle=0$).
\end{itemize}

Note that a hypersphere $S$ can also be interpreted as the set of
points $x\in S$. Correspondingly, it admits, along with the
representation $\hat{s}$, the dual representation as a transversal
intersection of $\bbP(\bbL^{N+1,1})$ with the projective $N$-space
$\bbP(\hat{s}^\perp)$, polar to the point $\hat{s}$ with respect
to $\bbP(\bbL^{N+1,1})$; here, of course,
$\hat{s}^\perp=\{\hat{x}\in\bbR^{N+1,1}: \langle
\hat{s},\hat{x}\rangle=0\}$. This can be generalized to model
lower-dimensional spheres.
\begin{itemize}
\item {\em Spheres.} A $k$-sphere is a (generic) intersection of
$N-k$ hyperspheres $S_i$ $\,(i=1,\ldots,N-k)$. The intersection of
$N-k$ hyperspheres represented by $\hat{s}_i\in\bbR^{N+1,1}_{\rm
out}$ $\,(i=1,\ldots,N-k)$ is generic if the $(N-k)$-dimensional
linear subspace of $\bbR^{N+1,1}$ spanned by $\hat{s}_i$ is
space-like:
\[
\Sigma={\rm
span}(\hat{s}_1,\ldots,\hat{s}_{N-k})\subset\bbR^{N+1,1}_{\rm
out}.
\]
As a set of points, this $k$-sphere is represented as
$\bbP(\bbL^{N+1,1}\cap\Sigma^\perp$, where
\[
\Sigma^\perp=\bigcap_{i=1}^{N-k}\hat{s}_i^\perp=\Big\{\hat{x}\in\bbR^{N+1,1}:\;
\langle\hat{s}_1,\hat{x}\rangle=\ldots=
\langle\hat{s}_{N-k},\hat{x}\rangle=0\Big\}
\]
is a $(k+2)$-dimensional linear subspace of $\bbR^{N+1,1}$ of
signature $(k+1,1)$.

Through any $k+2$ points $x_1,\ldots,x_{k+2}\in\bbR^N$ in general
position one can draw a unique $k$-sphere. It corresponds to the
$(k+2)$-dimensional linear subspace
\[
\Sigma^\perp={\rm span}(\hat{x}_1,\ldots,\hat{x}_{k+2}),
\]
of signature $(k+1,1)$, with $k+2$ linearly independent isotropic
vectors $\hat{x}_1,\ldots,\hat{x}_{k+2}\in\bbL^{N+1,1}$. In the
polar formulation, this $k$-sphere corresponds to the
$(N-k)$-dimensional space-like linear subspace
\[
\Sigma=\bigcap_{i=1}^{k+2}\hat{x}_i^\perp=\Big\{\hat{s}\in\bbR^{N+1,1}:\;
\langle\hat{s},\hat{x}_1\rangle=\ldots=
\langle\hat{s},\hat{x}_{k+2}\rangle=0\Big\}.
\]
\end{itemize}

M\"obius geometry is the study of properties of (non)-oriented
hyperspheres invariant with respect to projective transformations
of $\bbP(\bbR^{N+1,1})$ which map points to points, i.e., which
leave $\bbP(\bbL^{N+1,1})$ invariant. Such transformations are
called {\em M\"obius transformations}.
\begin{theorem}\label{thm: Moeb fund}
{\bf (Fundamental theorem of M\"obius geometry)}

\noindent a) The group of M\"obius transformations is isomorphic
to $O(N+1,1)/\{\pm I\}\simeq O^+(N+1,1)$, the group  of Lorentz
transformations of $\,\bbR^{N+1,1}$ preserving the time-like
direction.

\noindent b) Every conformal diffeomorphism of
$\,\bbS^N\simeq\bbR^N\cup\{\infty\}$ is induced by the restriction
to $\bbP(\bbL^{N+1,1})$ of a M\"obius transformation.
\end{theorem}

The group $O^+(N+1,1)$ is generated by reflections
\begin{equation}\label{eq: Moeb transf}
A_{\hat{s}}:\bbR^{N+1,1}\to\bbR^{N+1,1},\qquad
A_{\hat{s}}(\hat{x})=\hat{x}-\frac{2\langle\hat{s},\hat{x}\rangle}
{\langle\hat{s},\hat{s}\rangle}\,\hat{s}.
\end{equation}
If $\hat{s}$ is a hypersphere (\ref{eq: Moeb sphere}), then the
transformation induced on $\bbR^N$ by $A_{\hat{s}}$ is obtained
from (\ref{eq: Moeb transf}) by a computation with the
representatives (\ref{eq: Moeb point}) for points and is given by:
\begin{equation}\label{eq: inversion}
x\mapsto c+\frac{r^2}{|x-c|^2}\,(x-c)
\end{equation}
(inversion in the hypersphere $S=\{x\in\bbR^N:|x-c|^2=r^2\}$);
similarly, if $\hat{s}=\hat{p}\,$ is the hyperplane (\ref{eq: Moeb
plane}), then the transformation induced on $\bbR^N$ by
$A_{\hat{p}}$ is easily computed to be
\begin{equation}\label{eq: reflection}
x\mapsto x-\frac{2(\langle v,x\rangle-d)}{\langle v,v\rangle}\,v
\end{equation}
(reflection in the hyperplane $P=\{x\in\bbR^N:\langle
v,x\rangle=d\}$).

Since (non-)vanishing of the $\ee_\infty$-component of a point in
$\bbP(\bbR^{N+1,1})$ is not invariant under a general M\"obius
transformation, there is no distinction in M\"obius geometry
between hyperspheres and hyperplanes.

\subsection{Laguerre geometry}
\label{subsect: Laguerre}

Also in the case of {\em Laguerre geometry} the Blaschke's book
\cite{BlIII} serves as the indispensable classical source. One can
find a modern account, e.g., in \cite{Benz, Ce, PP}.

Laguerre geometry is a subgeometry of Lie geometry, with
hyperplanes distinguishable among all hyperspheres, as
hyperspheres through $\infty$. Thus, Laguerre geometry studies
properties of hyperspheres invariant under the subgroup of Lie
sphere transformations which preserve the set of hyperplanes.
Following objects in $\bbR^N$ are elements of the Laguerre
geometry.

\begin{itemize}
\item {\em (Oriented) hyperspheres} $S=\{x\in\bbR^N:
|x-c|^2=r^2\}$ with centers $c\in\bbR^N$ and signed radii
$r\in\bbR$, can be put into correspondence with $(N+1)$-tuples
$(c,r)$.

\item{\em Points} $x\in\bbR^N$ are considered as hyperspheres of
radius zero, and are put into correspondence with $(N+1)$-tuples
$(x,0)$.

\item {\em (Oriented) hyperplanes} $P=\{x\in\bbR^N: \langle
v,x\rangle=d\}$, with unit normals $v\in\bbS^{N-1}$ and
$d\in\bbR$, can be put into correspondence with $(N+1)$-tuples
$(v,d)$.
\end{itemize}

In the projective model of Lie geometry, hyperplanes are
distinguished as elements of $\bbP(\bbL^{N+1,2})$ with the
vanishing $\ee_{0}$-component. (Of course, one could replace here
$\ee_{0}$ by any isotropic vector.) Thus, Laguerre geometry
studies the subgroup of Lie sphere transformations preserving the
subset of $\bbP(\bbL^{N+1,2})$ with the vanishing
$\ee_{0}$-component.

There seems to exist no model of Laguerre geometry, where
hyperspheres and hyperplanes would be modelled as points of one
and the same space. Depending on which of both types of elements
is modelled by points, one comes to the {\em Blaschke cylinder
model} or to the {\em cyclographic model} of Laguerre geometry. We
will use the first of these models which has an advantage of a
simpler description of the distinguished objects of the Laguerre
geometry, which are hyperplanes. The main advantage of the second
model is a simpler description of the group of Laguerre
transformations.
\smallskip

The scene of the both models consists of {\em two}
$(N+1)$-dimensional projective spaces, whose spaces of homogeneous
coordinates, $\bbR^{N,1,1}$ and $(\bbR^{N,1,1})^*$, are dual to
one another and arise from $\bbR^{N+1,2}$ by ``forgetting'' the
$\ee_0$-, resp. $\ee_\infty$-components. Thus, $\bbR^{N,1,1}$ is
spanned by $N+2$ linearly independent vectors
$\ee_1,\ldots,\ee_N$, $\ee_{N+3}$, $\ee_\infty$, and is equipped
with a degenerate bilinear form of the signature $(N,1,1)$ in
which the above vectors are pairwise orthogonal, the first $N$
being space-like: $\langle \ee_i,\ee_i\rangle=1$ for $1\le i\le
N$, while the last two being time-like and isotropic,
respectively: $\langle \ee_{N+3},\ee_{N+3}\rangle=-1$ and $\langle
\ee_\infty,\ee_\infty\rangle=0$. Similarly, $(\bbR^{N,1,1})^*$ is
assumed to have an orthogonal basis consisting of
$\ee_1,\ldots,\ee_N$, $\ee_{N+3}$, $\ee_0$, again with an
isotropic last vector: $\langle \ee_0,\ee_0\rangle=0$. Note that
one and the same symbol $\langle\cdot,\cdot\rangle$ is used to
denote two degenerate bilinear forms in our two spaces. We will
overload this symbol even more and use it also for the
(non-degenerate) pairing between these two spaces, which is
established by setting $\langle
\ee_0,\ee_\infty\rangle=-\tfrac{1}{2}$, additionally to the above
relations. (Note that a degenerate bilinear form cannot be used to
identify a vector space with its dual.)

In both models mentioned above there holds:
\begin{itemize}
\item {\em Hyperplane $P=(v,d)$ is modelled as a point in the
space $\bbP(\bbR^{N,1,1})$ with a representative}
\begin{equation}\label{eq: Bl plane}
\hat{p}=v+2d\ee_\infty+\ee_{N+3}.
\end{equation}
\item {\em Hypersphere $S=(c,r)$ is modelled as a point in the
space $\bbP\big((\bbR^{N,1,1})^*\big)$ with a representative}
\begin{equation}\label{eq: Lag sphere}
\hat{s}=c+\ee_0+r\ee_{N+3}.
\end{equation}
\end{itemize}
Each one of the models appears if one considers one of the spaces
as a preferred (fundamental) one, and interprets the points of the
second space as hyperplanes in the preferred one. In the {\em
Blaschke cylinder model}, the preferred space is the space
$\bbP(\bbR^{N,1,1})$ whose points model hyperplanes
$P\subset\bbR^N$. A hypersphere $S\subset\bbR^N$ is then modelled
as a hyperplane $\{\xi\in\bbP(\bbR^{N,1,1}): \langle
\hat{s},\xi\rangle=0\}$ in the space $\bbP(\bbR^{N,1,1})$. Basic
features of this model:
\begin{itemize}
\item[(i)] Oriented hyperplanes $P\subset\bbR^N$ are in a
one-to-one correspondence with points $\hat{p}$ of the quadric
$\bbP(\bbL^{N,1,1})$, where
\begin{equation}\label{eq: Lag L}
\bbL^{N,1,1}=\big\{\xi\in\bbR^{N,1,1}:\langle
\xi,\xi\rangle=0\big\}.
\end{equation}

\item[(iii)] Two oriented hyperplanes $P_1,P_2\subset\bbR^N$ are
in an oriented contact (parallel), if and only if their
representatives $\hat{p}_1$, $\hat{p}_2$ differ by a vector
parallel to $\ee_\infty$.

\item[(iii)] An oriented hypersphere $S\subset\bbR^N$ is in an
oriented contact with an oriented hyperplane $P\subset\bbR^N$, if
and only if $\hat{p}\in\hat{s}$, that is, if $\langle
\hat{p},\hat{s}\rangle=0$. Thus, a hypersphere $S$ is interpreted
as a set of all its tangent hyperplanes.
\end{itemize}
The quadric $\bbP(\bbL^{N,1,1})$ is diffeomorphic to the {\em
Blaschke cylinder}
\begin{equation}\label{eq: Bl cyl}
\cZ=\big\{(v,d)\in\bbR^{N+1}:
|v|=1\big\}=\bbS^{N-1}\times\bbR\subset\bbR^{N+1}.
\end{equation}
Two points of this cylinder represent parallel hyperplanes, if
they lie on one straight line generator of $\cZ$ parallel to its
axis. In the ambient space $\bbR^{N+1}$ of the Blaschke cylinder,
oriented hyperspheres $S\subset\bbR^N$ are in a one-to-one
correspondence with hyperplanes non-parallel to the axis of $\cZ$:
\begin{equation}
S\sim\big\{(v,d)\in\bbR^{N+1}:\langle c,v\rangle-d-r=0\big\}.
\end{equation}
An intersection of such a hyperplane with $\cZ$ consists of points
in $\cZ$ which represent tangent hyperplanes to $S\subset\bbR^N$,
as follows from eq. (\ref{eq: tangency}).

In this paper, we will not use the cyclographic model of Laguerre
geometry; its short description is put in Appendix \ref{Subsect:
Laguerre cyclographic}.

\section{Discrete curvature line parametrization in Lie, M\"obius
and Laguerre geometries} \label{Sect: curv}

Starting from here, we restrict ourselves to the geometry of
surfaces in the three-dimensional Euclidean space $\bbR^3$.
Accordingly, one should set $N=3$ in all previous considerations.

It is natural to consider following objects as discrete surfaces
in the various geometries discussed above:
\begin{itemize}
\item In {\em Lie geometry}, a surface is viewed as built of its
contact elements. These contact elements are interpreted as points
of the surface and tangent planes (or, equivalently, normals) at
these points. This can be discretized in a natural way: a discrete
surface is a map
\[
(x,P):\,\bbZ^2\to\{\rm contact\ elements\ of\ surfaces\ in\
\bbR^3\},
\]
or, in the projective model of Lie geometry, a map
\begin{equation}\label{eq: surf Lie}
\ell:\,\bbZ^2\to\cL_0^{4,2},
\end{equation}
where, recall, $\cL_0^{4,2}$ denotes the set of isotropic lines in
$\bbP(\bbR^{4,2})$.

\item In {\em M\"obius geometry}, a surface is viewed simply as
built of points. A discrete surface is a map
\[
x:\,\bbZ^2\to\bbR^3,
\]
or, in the projective model, a map
\begin{equation}\label{eq: surf Moeb}
\hat{x}:\,\bbZ^2\to\bbP(\bbL^{4,1}).
\end{equation}

\item In {\em Laguerre geometry}, a surface is viewed as the
envelope of the system of its tangent planes. A discrete surface
is a map
\[
P:\,\bbZ^2\to\{\rm oriented\ planes\ in\ \bbR^3\},
\]
or, in the projective model, a map
\begin{equation}\label{eq: surf Lag}
\hat{p}:\,\bbZ^2\to\bbP(\bbL^{3,1,1}).
\end{equation}
It should be mentioned that a substantial part of the description
of a surface in Laguerre geometry is its {\em Gauss map}
\begin{equation}\label{eq: Gauss map}
v:\,\bbZ^2\to\bbS^2,
\end{equation}
consisting of unit normals $v$ to the tangent planes $P=(v,d)$.
\end{itemize}

Thus, description of a discrete surface in Lie geometry contains
more information than description of a discrete surface in
M\"obius or in Laguerre geometry. Actually, the former merges the
two latter ones.

\subsection{Lie geometry}
\label{subsect: Lie curv}

The following definition is a discretization of the Lie-geometric
description of curvature line parametrized surfaces, as found,
e.g., in \cite{BlIII}.

\begin{definition}\label{def: Lie curv Euc}
{\bf (Principal contact element nets. Euclidean model)} A map
\[
(x,P):\bbZ^2\to\{\rm contact\ elements\ of\ surfaces\ in\ \bbR^3\}
\]
is called a principal contact element net, if any two neighboring
contact elements $(x,P)$, $(x_i,P_i)$ have a sphere $S^{(i)}$ in
common, that is, a sphere touching both planes $P$, $P_i$ at the
corresponding points $x$, $x_i$.
\end{definition}
Thus, the normals to the neighboring planes $P$, $P_i$ at the
corresponding points $x$, $x_i$ intersect at a point $c^{(i)}$
(the center of the sphere $S^{(i)}$), and the distances from
$c^{(i)}$ to $x$ and to $x_i$ are equal, see Fig.~\ref{Fig: curv
sph}. The spheres $S^{(i)}$, attached to the edges of $\bbZ^2$
parallel to the $i$-th coordinate axis, will be called {\em
principal curvature spheres} of the discrete surface.

\begin{figure}[htbp]
 \psfrag{x}[Bl][bl][0.9]{$x$}
 \psfrag{P}[Bl][bl][0.9]{$P$}
 \psfrag{x_i}[Bl][bl][0.9]{$x_i$}
 \psfrag{P_i}[Bl][bl][0.9]{$P_i$}
 \psfrag{S^i}[Bl][bl][0.9]{$S^{(i)}$}
 \psfrag{c^i}[Bl][bl][0.9]{$c^{(i)}$}
 \center{\includegraphics[width=120mm]{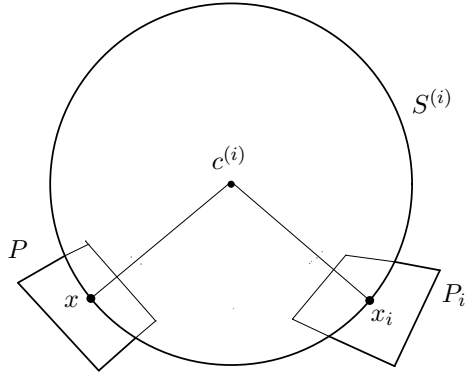}}
 \caption{Principal curvature sphere}
 \label{Fig: curv sph}
\end{figure}

A direct translation of Definition \ref{def: Lie curv Euc} into
the projective model looks as follows:
\begin{definition}\label{def: Lie curv}
{\bf (Principal contact element nets. Projective model)} A map
$\ell:\bbZ^2\to\cL_0^{4,2}$ is called a principal contact element
net, if it is a discrete congruence of isotropic lines in
$\bbP(\bbR^{4,2})$, that is, any two neighboring lines intersect:
\begin{equation}\label{eq: Lie focal}
\ell(u)\cap\ell(u+e_i)=\hat{s}^{(i)}(u)\in\bbP(\bbL^{4,2}), \quad
\forall u\in\bbZ^2,\quad \forall i=1,2.
\end{equation}
\end{definition}

In the projective model, the representatives of the principal
curvature spheres $S^{(i)}$ of the $i$-th coordinate direction
build the corresponding focal net of the line congruence $\ell$;
\begin{equation}\label{eq: Lie focal nets}
\hat{s}^{(i)}:\bbZ^2\to\bbP(\bbL^{4,2}),\quad i=1,2,
\end{equation}
cf. Definition \ref{Def: focal net}. According to Theorem
\ref{thm: focal Q-nets}, both focal nets are Q-nets in
$\bbP(\bbR^{4,2})$. This motivates the following definition.
\begin{definition}\label{def: R-congr}
{\bf (Discrete R-congruence of spheres)} A map
\[
S:\bbZ^m\to\{\rm oriented\ spheres\ in\ \bbR^3\}
\]
is called a discrete R-congruence (Ribaucour congruence) of
spheres, if the corresponding map
\[
\hat{s}:\bbZ^m\to\bbP(\bbL^{4,2})
\]
is a Q-net in $\bbP(\bbR^{4,2})$.
\end{definition}
A geometric characterization of discrete R-congruences will be
given in Sect. \ref{sect: R-congr}.
\begin{corollary}\label{cor: curv spheres}
{\bf (Curvature spheres build an R-congruence)} For a discrete
contact element net, the principal curvature spheres of the $i$-th
coordinate direction (i=1,2) build a two-dimensional discrete
R-congruence.
\end{corollary}

Turning to transformations of principal contact element nets, we
introduce the following definition.
\begin{definition}\label{def: Lie Ribaucour Euc}
{\bf (Ribaucour transformation. Euclidean model)} Two principal
contact element nets
\[
(x,P), (x^+\!,P^+):\bbZ^2\to\{\rm contact\ elements\ of\ surfaces\
in\ \bbR^3\}
\]
are called Ribaucour transforms of one another, if any two
corresponding contact elements $(x,P)$ and $(x^+\!,P^+)$ have a
sphere $S$ in common, that is, a sphere which touches both planes
$P$, $P^+$ at the corresponding points $x$, $x^+$.
\end{definition}

\begin{figure}[htbp]
 \psfrag{x}[Bl][bl][0.9]{$x$}
 \psfrag{P}[Bl][bl][0.9]{$P$}
 \psfrag{x^+}[Bl][bl][0.9]{$x^+$}
 \psfrag{P^+}[Bl][bl][0.9]{$P^+$}
 \psfrag{S}[Bl][bl][0.9]{$S$}
 \center{\includegraphics[width=120mm]{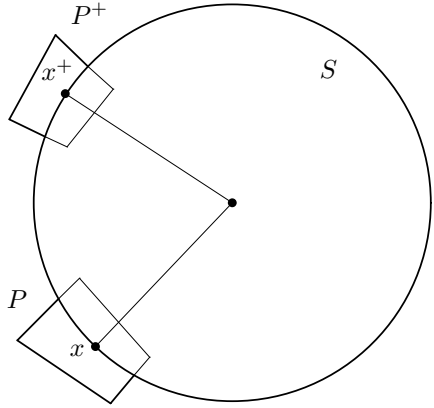}}
 \caption{Ribaucour transformation}
 \label{Fig: Rib sph}
\end{figure}

Again, a direct translation of Definition \ref{def: Lie Ribaucour
Euc} into the language of the projective model gives:

\begin{definition}\label{def: Lie Ribaucour}
{\bf (Ribaucour transformation. Projective model)} Two principal
contact element nets $\ell,\ell^+:\bbZ^2\to\cL_0^{4,2}$ are called
Ribaucour transforms of one another, if these discrete congruences
of isotropic lines are in the relation of F-transformation, that
is, if any pair of the corresponding lines intersect:
\begin{equation}\label{eq: Lie Ribaucour}
\ell(u)\cap\ell^+(u)=\hat{s}(u)\in\bbP(\bbL^{4,2}),\quad \forall
u\in\bbZ^2.
\end{equation}
\end{definition}

Spheres $S$ of a Ribaucour transformation are attached to the
vertices $u$ of the lattice $\bbZ^2$, or, better, to the
``vertical'' edges connecting the vertices $(u,0)$ and $(u,1)$ of
the lattice $\bbZ^2\times\{0,1\}$. In the projective model, their
representatives
\begin{equation}
\hat{s}:\bbZ^2\to\bbP(\bbL^{4,2})
\end{equation}
build the focal net of the three-dimensional line congruence for
the third coordinate direction. From Theorem \ref{thm: focal
Q-nets} there follows:

\begin{corollary}\label{cor: Ribaucour spheres}
{\bf (Spheres of a Ribaucour transformation build an
R-congruence)} The spheres of a generic Ribaucour transformation
build a discrete R-congruence.
\end{corollary}

Now, we turn to the study of the geometry of an elementary
quadrilateral of contact elements of a principal contact element
net, consisting of $\ell\sim (x,P)$, $\ell_1\sim (x_1,P_1)$,
$\ell_2\sim (x_2,P_2)$, and $\ell_{12}\sim (x_{12}, P_{12})$.

We leave aside a degenerate {\em umbilic} situation, when all four
lines have a common point and span a four-dimensional space.
Geometrically, this means that one is dealing with four contact
elements of a sphere $S\subset\bbR^3$. In this situation, one
cannot draw any further conclusion about the four points
$x,x_1,x_2,x_{12}$ on the sphere $S$: they can be arbitrary.

In the non-umbilic situation, the space spanned by the four lines
$\ell,\ell_1,\ell_2,\ell_{12}$ is three-dimensional. The four
elements
$\hat{x},\hat{x}_1,\hat{x}_2,\hat{x}_{12}\in\bbP(\bbL^{4,2})$
corresponding to the points $x,x_1,x_2,x_{12}\in\bbR^3$ are
obtained as the intersection of the four isotropic lines
$\ell,\ell_1,\ell_2,\ell_{12}$ with the projective hyperplane
$\bbP(\ee_6^\perp)$ in $\bbP(\bbR^{4,2})$. Therefore, the four
elements $\hat{x},\hat{x}_1,\hat{x}_2,\hat{x}_{12}$ lie in a
plane. A suitable framework for the study of this configuration is
the projective model of the M\"obius geometry. Namely, omitting
the inessential (vanishing) $\ee_6$-component, we arrive at a
planar quadrilateral in the M\"obius sphere $\bbP(\bbL^{4,1})$. We
devote Sect. \ref{subsect: O-nets} to the study of such objects.

Analogously, the four elements
$\hat{p},\hat{p}_1,\hat{p}_2,\hat{p}_{12}\in\bbP(\bbL^{4,2})$
corresponding to the planes $P,P_1,P_2,P_{12}\in\bbR^3$ are
obtained as the intersection of the four isotropic lines
$\ell,\ell_1,\ell_2,\ell_{12}$ with the projective hyperplane
$\bbP(\ee_\infty^\perp)$ in $\bbP(\bbR^{4,2})$. Therefore, also
the four elements $\hat{p},\hat{p}_1,\hat{p}_2,\hat{p}_{12}$ lie
in a plane. A suitable framework for the study of such a
configuration is the projective model of the Laguerre geometry;
this will be performed in Sect. \ref{subsect: conical}.

\subsection{M\"obius geometry: circular nets}
\label{subsect: O-nets}

Circular nets were introduced and studied in the context of
integrable systems in \cite{B, CDS, KSch}.
\smallskip

{\em Caution:} in this section, the notation $\hat{x}$ refers to
the M\"obius-geometric representatives in $\bbL^{4,1}$, and not to
the Lie-geometric ones in $\bbL^{4,2}$. The former are obtained
from the latter one by omitting the (vanishing) $\ee_6$-component.
\smallskip

We assume that the principal contact element nets under
consideration are generic, i.e., do not contain umbilic
quadruples. The main result of this section is the following
claim.
\begin{theorem} \label{thm: dos}
{\bf (Points of principal contact element nets form circular
nets)} For a principal contact element net
\[
(x,P):\bbZ^2\to\{\rm contact\ elements\ of\ surfaces\ in\
\bbR^3\},
\]
its points $x:\bbZ^2\to\bbR^3$ form a circular net.
\end{theorem}

This statement refers to the notion, which can be defined in two
different ways.
\begin{definition}\label{def: O-net}
{\bf (Circular net. Euclidean model)} A net $x:\bbZ^m\to\bbR^3$ is
called circular, if the vertices of any elementary quadrilateral
$(x,x_i,x_{ij},x_j)$ (at any $u\in\bbZ^m$ and for all pairs $1\le
i\neq j\le m$) lie on a circle (in particular, are co-planar).
\end{definition}

\begin{definition}\label{def: O-net in Moeb}
{\bf (Circular net. Projective model)} A net $x:\bbZ^m\to\bbR^3$
is called circular, if the corresponding
$\hat{x}:\bbZ^m\to\bbP(\bbL^{4,1})$ is a Q-net in
$\bbP(\bbR^{4,1})$.
\end{definition}

This time a translation between the Euclidean model and the
projective model is not straightforward and actually constitutes
the matter of Theorem \ref{thm: dos}: indeed, this theorem has
already been demonstrated (or, better, is obvious) in terms of
Definition \ref{def: O-net in Moeb}, and it remains to establish
the equivalence of Definitions \ref{def: O-net},\ref{def: O-net in
Moeb}.
\smallskip

\noindent {\bf Conceptual proof.} The linear subspace of
$\bbR^{4,1}$ spanned by the isotropic vectors $\hat{x}$,
$\hat{x}_i$, $\hat{x}_j$, $\hat{x}_{ij}$ is three-dimensional. Its
orthogonal complement is therefore two-dimensional and lies in
$\bbR^{4,1}_{\rm out}$. Therefore, it represents a circle (an
intersection of two spheres). $\Box$
\smallskip

\noindent {\bf Computational proof.} For arbitrary representatives
$\tilde{x}\in\bbL^{4,1}$ of $\hat{x}$, the requirement of
Definition \ref{def: O-net in Moeb} is equivalent to equation of
the type (\ref{eq: dcn hom}). Since the representatives
$\hat{x}=x+\ee_0+|x|^2\ee_\infty$ fixed in (\ref{eq: Moeb point})
lie in an affine hyperplane of $\bbR^{4,1}$ (their
$\ee_0$-component is equal to 1), one has an equation of the type
(\ref{eq: dcn property}) for them. Clearly, this holds if and only
if $x$ is a Q-net in $\bbR^3$ and $|x|^2$ satisfies the same
equation (\ref{eq: dcn property}) as $x$ does. We show that the
latter condition is equivalent to circularity. On a single planar
elementary quadrilateral $(x,x_i,x_{ij},x_j)$, the function
$|x|^2$ satisfies eq. (\ref{eq: dcn property}) simultaneously with
$|x-c|^2=|x|^2-2\langle x,c\rangle+|c|^2$ with any $c\in\bbR^3$.
Choose $c$ to be the center of the circle through the three points
$x$, $x_i$, $x_j$, so that $|x-c|^2=|x_i-c|^2=|x_j-c|^2$. Then eq.
(\ref{eq: dcn property}) for $|x-c|^2$ turns into
$|x_{ij}-c|^2=|x-c|^2$, which means that $x_{ij}$ lies on the same
circle. $\Box$
\medskip

\begin{figure}[htbp]
\begin{center}
\epsfig{file=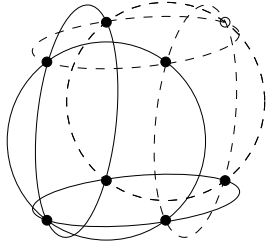,width=50mm} \caption{An elementary
hexahedron of a circular net} \label{Fig: Miquel 3d}
\end{center}
\end{figure}

Two-dimensional circular nets ($m=2$) are discrete analogs of the
curvature lines parametrized surfaces, while the case $m=3$
discretizes orthogonal coordinate systems in $\bbR^3$. A
construction of an elementary hexahedron of a circular net is
based on the following geometric theorem:
\begin{theorem} \label{Th: Miquel 3d}
{\bf (Elementary hexahedron of a circular net)} Given seven points
$x$, $x_i$, and $x_{ij}$ $\,(1\le i<j\le 3)$ in $\,\bbR^3$, such
that each of the three quadruples $(x,x_i,x_j,x_{ij})$ lies on a
circle $C_{ij}$, define three new circles $\tau_iC_{jk}$ as those
passing through the triples $(x_i,x_{ij},x_{ik})$, respectively.
Then these new circles intersect at one point, see Fig. \ref{Fig:
Miquel 3d}:
\[
x_{123}=\tau_1C_{23}\cap\tau_2C_{31}\cap\tau_3C_{12}\,.
\]
\end{theorem}

{\bf Proof.} This is a particular case of Theorem \ref{Th: Q-net
in Q}, applied to the quadric $\bbP(\bbL^{4,1})$. $\Box$
\smallskip

This theorem can be proven also by elementary geometrical
considerations. If one notes that under conditions of Theorem
\ref{Th: Miquel 3d} the seven points $x$, $x_i$, $x_{ij}\,$ lie on
a two-dimensional sphere, and performs a stereographic projection
of this sphere with the pole at $x$, one arrives at a planar
picture which is nothing but the classical Miquel theorem.

\subsection{Laguerre geometry: conical nets}
\label{subsect: conical}

Conical meshes have been introduced recently in \cite{LPWYW}.
\smallskip

{\em Caution:} in this section, the notation $\hat{p}$ refers to
the Laguerre-geometric representatives in $\bbL^{3,1,1}$, and not
to the Lie-geometric ones in $\bbL^{4,2}$. The former are obtained
from the latter by omitting the (vanishing) $\ee_0$-component.
\smallskip

As in the previous section, we assume that the principal contact
element nets under consideration do not contain umbilic
quadruples. The main result of this section is the following
claim.
\begin{theorem} \label{thm: conical}
{\bf (Tangent planes of principal contact element nets form
conical nets)} For a principal contact element net
\[
(x,P):\bbZ^2\to\{\rm contact\ elements\ of\ surfaces\ in\
\bbR^3\},
\]
its tangent planes $P:\bbZ^m\to\{\rm oriented\ planes\ in\
\bbR^3\}$ form a conical net.
\end{theorem}

This statement refers to the notion which can be defined in two
different ways.
\begin{definition}\label{def: conical}
{\bf (Conical net. Euclidean model)} A net $P:\bbZ^m\to\{\rm
oriented\ planes \ in \ \bbR^3\}$ is called conical, if at any
$u\in\bbZ^m$ and for all pairs $1\le i\neq j\le m$ the four planes
$P,P_i,P_{ij},P_j$ touch a cone of revolution (in particular,
intersect at the tip of the cone).
\end{definition}

\begin{definition}\label{def: conical in Lag}
{\bf (Conical net. Projective model)} A net $P:\bbZ^m\to\{\rm
oriented\ planes \ in\ \bbR^3\}$ is called conical, if the
corresponding $\hat{p}:\bbZ^m\to\bbP(\bbL^{3,1,1})$ is a Q-net in
$\bbP(\bbR^{3,1,1})$.
\end{definition}

Theorem \ref{thm: conical} is obvious in terms of Definition
\ref{def: conical in Lag}, so the real content of this theorem is
the translation between the Euclidean model and the projective
model, that is, establishing equivalence of Definitions \ref{def:
conical}, \ref{def: conical in Lag}.
\smallskip

{\bf Proof.} Representatives $\hat{p}$ in (\ref{eq: Bl plane})
form a Q-net, if and only if they satisfy eq. (\ref{eq: dcn
property}), that is, if $v:\bbZ^m\to\bbS^2$ and $d:\bbZ^m\to\bbR$
satisfy this equation. Equation (\ref{eq: dcn property}) for $v$
yields that $v:\bbZ^m\to\bbS^2$ is actually a Q-net in $\bbS^2$,
so that any quadrilateral $(v,v_i,v_{ij},v_j)$ in $\bbS^2$ is
planar and therefore circular. Equation (\ref{eq: dcn property})
for $(v,d)$ yields that the (unique) intersection point of three
planes $P$, $P_i$, $P_j$ lies on $P_{ij}$, as well, so that all
four planes intersect in one point. Thus, we arrived at a
characterization of conical nets in the sense of Definition
\ref{def: conical in Lag} as those nets of planes for which every
quadruple of planes $(P,P_i,P_{ij},P_j)$ is concurrent and every
quadrilateral $(v,v_i,v_{ij},v_j)$ of unit normal vectors is
planar. It is clear that this description is equivalent to that of
Definition \ref{def: conical}. The direction of the axis of the
tangent cone coincides with the spherical center of the
quadrilateral $(v,v_i,v_{ij},v_j)$ in $\bbS^2$. $\Box$
\smallskip

Thus, conical nets are Q-nets with circular Gauss maps. It is
worthwhile to mention that, in order to prescribe a conical net,
it is enough to prescribe a circular Gauss map $v:\bbZ^m\to\bbS^2$
and additionally the numbers $d$ (i.e., the planes $P=(v,d)$)
along the coordinate axes of $\bbZ^m$. Indeed, these data allow
one to reconstruct the conical net uniquely. This is done via a
recursive procedure, whose elementary step consists in finding the
fourth plane $P_{ij}$ provided three planes $P,P_i,P_j$ and the
normal direction $v_{ij}$ of the fourth one are known. But this is
easy: $P_{ij}$ is the plane normal to $v_{ij}$ through the unique
intersection point of the three planes $P,P_i,P_j$.

\subsection{Synthesis}
\label{subsect: synthesis}

In view of Theorems \ref{thm: dos}, \ref{thm: conical}, it is
natural to ask whether, given a circular net $x:\bbZ^2\to\bbR^3$,
or a conical net $P:\bbZ^2\to\{\rm oriented\ planes\ in\
\bbR^3\}$, there exists a principal contact element net
\[
(x,P):\bbZ^2\to\{\rm contact\ elements\ of\ surfaces\ in\
\bbR^3\},
\]
with the prescribed half of the data ($x$ or $P$). A positive
answer to this question is a corollary of the following general
theorem.

\begin{theorem}\label{thm: synthesis}
{\bf (Extending R-congruences of spheres to curvature line
parametrized surfaces)} Given a discrete R-congruence of spheres
\[
S:\bbZ^2\to\{{\rm oriented\ spheres\ in\ }\bbR^3\},
\]
there exists a two-parameter family of principal contact element
nets
\[
(x,P):\bbZ^2\to\{\rm contact\ elements\ of\ surfaces\ in\ \bbR^3\}
\]
such that $S$ belongs to the contact element $(x,P)$, i.e., $P$ is
the tangent plane to $S$ at the point $x\in S$, for all
$u\in\bbZ^2$. Such a principal contact element net is uniquely
determined by prescribing a contact element $(x,P)(0,0)$
containing the sphere $S(0,0)$.
\end{theorem}
{\bf Proof.} The input data is a Q-net
$\hat{s}:\bbZ^2\to\bbP(\bbL^{4,2})$ in the Lie quadric, and we are
looking for a congruence of isotropic lines
$\ell:\bbZ^2\to\cL_0^{4,2}$ such that $\hat{s}(u)\in\ell(u)$ for
all $u\in\bbZ^2$. The construction starts with an arbitrary
isotropic line $\ell(0,0)$ through $\hat{s}(0,0)$, and hinges on
the following lemma.
\begin{lemma}\label{lemma: synthesis}
For an isotropic line $\ell\in\cL_0^{4,2}$ and a point
$\hat{s}_1\in\bbP(\bbL^{4,2})$ not lying on $\ell$, there is a
unique isotropic line $\ell_1$ through $\hat{s}_1$ intersecting
$\ell$.
\end{lemma}
{\bf Proof.} Let $\hat{s},\hat{\sigma}$ be two arbitrary points on
$\ell$ (in homogeneous coordinates), so that the line $\ell$ is
given by the linear combinations
$\alpha\hat{s}+\beta\hat{\sigma}$. Relation
$\langle\alpha\hat{s}+\beta\hat{\sigma},\hat{s}_1\rangle=0$ yields
\[
\alpha:\beta=-\,\langle\hat{\sigma},\hat{s}_1\rangle:\langle
\hat{s},\hat{s}_1\rangle.
\]
Thus, there exists a unique point $\hat{s}^{(1)}\in\ell$ such that
$\langle \hat{s}^{(1)},\hat{s}_1\rangle=0$. Now $\ell_1$ is the
line through $\hat{s}_1$ and $\hat{s}^{(1)}$. $\Box$
\smallskip

\noindent {\bf Proof of Theorem \ref{thm: synthesis}, continued.}
With the help of Lemma \ref{lemma: synthesis}, one can construct
the isotropic lines of the congruence along the coordinate axes,
\[
\ell:\bbZ\times\{0\}\to\cL_0^{4,2}\quad {\rm and}\quad
\ell:\{0\}\times\bbZ\to\cL_0^{4,2}.
\]
Next, one has to extend the congruence $\ell$ from the coordinate
axes to the whole of $\bbZ^2$. An elementary step of this
extension consists in finding, for three given isotropic lines
$\ell$, $\ell_1$, $\ell_2$ (such that $\ell$ intersects both
$\ell_1$ and $\ell_2$) the fourth one, $\ell_{12}$, intersecting
$\ell_1$ and $\ell_2$ and going through a given point
$\hat{s}_{12}$. One can use for this Lemma \ref{lemma: synthesis},
but then one has to demonstrate that this construction is
consistent, i.e., that the lines $\ell_{12}$ obtained from the
requirements of intersecting with $\ell_1$ and with $\ell_2$
coincide. We show this with the following argument. The space
$V={\rm span}(\ell,\ell_1,\ell_2)$ is three-dimensional. The
points $\hat{s},\hat{s}_1,\hat{s}_2$ lie in $V$. By the hypothesis
of the theorem, the quadrilateral
$(\hat{s},\hat{s}_1,\hat{s}_{12},\hat{s}_2)$ is planar, therefore
$\hat{s}_{12}$ lies in $V$, as well. Draw two planes in $V$:
$\Pi_1={\rm span}(\ell_1,\hat{s}_{12})$ and $\Pi_2={\rm
span}(\ell_2,\hat{s}_{12})$. Their intersection is a line
$\ell_{12}$ through $\hat{s}_{12}$. It remains to prove that this
line is isotropic. For this, note that $\ell_{12}$ can be
alternatively described as the line through two points
$\hat{s}_1^{(2)}=\ell_1\cap\ell_{12}$ and
$\hat{s}_2^{(1)}=\ell_2\cap\ell_{12}$. Both these points lie in
$\bbP(\bbL^{4,2})$, since they belong to the isotropic lines
$\ell_1$ and $\ell_2$, respectively. But it is easy to see that a
line in $\bbP(\bbR^{4,2})$ through two points from
$\bbP(\bbL^{4,2})$ is either isotropic, or contains no further
points from $\bbP(\bbL^{4,2})$, depending on whether these two
points are polar to one another (with respect to
$\bbP(\bbL^{4,2})$) or not. In our case the line $\ell_{12}$
contains, by construction, one further point $\hat{s}_{12}$ from
$\bbP(\bbL^{4,2})$, therefore it has to be isotropic. $\Box$
\medskip

Since the representatives $\hat{x}$ in $\bbP(\bbL^{4,2})$ of a
circular net $x:\bbZ^2\to\bbR^3$ form a Q-net in
$\bbP(\bbR^{4,2})$, and the same holds for the representatives
$\hat{p}$ in $\bbP(\bbL^{4,2})$ of a conical net
$P:\bbZ^2\to\{{\rm oriented\ planes\ in\ }\bbR^3\}$, we come to
the following conclusion (obtained independently by Pottmann
\cite{P}).
\begin{corollary}\label{cor: synthesis}
{\bf (Extending circular and conical nets to principal contact
element nets)}

{\rm i)} Given a circular net $x:\bbZ^2\to\bbR^3$, there exists a
two-parameter family of conical nets $P:\bbZ^2\to\{\rm planes\ in\
\bbR^3\}$ such that $x\in P$ for all $u\in\bbZ^2$, and the contact
element net
\[
(x,P):\bbZ^2\to\{\rm contact\ elements\ of\ surfaces\ in\ \bbR^3\}
\]
is principal. Such a conical net is uniquely determined by
prescribing a plane $P(0,0)$ through the point $x(0,0)$.

{\rm ii)} Given a conical net $P:\bbZ^2\to\{\rm oriented\ planes\
in\ \bbR^3\}$, there exists a two-parameter family of circular
nets $x:\bbZ^2\to\bbR^3$ such that $x\in P$ for all $u\in\bbZ^2$,
and the contact element net
\[
(x,P):\bbZ^2\to\{\rm contact\ elements\ of\ surfaces\ in\ \bbR^3\}
\]
is principal. Such a circular net is uniquely determined by
prescribing a point $x(0,0)$ in the plane $P(0,0)$.
\end{corollary}

\begin{figure}[htbp]
 \psfrag{s^1}[Bl][bl][0.9]{$S^{(1)}$}
 \psfrag{s_1^2}[Bl][bl][0.9]{$S_1^{(2)}$}
 \psfrag{s^2}[Bl][bl][0.9]{$S^{(2)}$}
 \psfrag{s_2^1}[Bl][bl][0.9]{$S_2^{(1)}$}
 \psfrag{l}[Bl][bl][0.9]{$\ell$}
 \psfrag{l_1}[Bl][bl][0.9]{$\ell_1$}
 \psfrag{l_2}[Bl][bl][0.9]{$\ell_2$}
 \psfrag{l_12}[Bl][bl][0.9]{$\ell_{12}$}
 \psfrag{x}[Bl][bl][0.9]{$x$}
 \psfrag{x_1}[Bl][bl][0.9]{$x_1$}
 \psfrag{x_2}[Bl][bl][0.9]{$x_2$}
 \psfrag{x_12}[Bl][bl][0.9]{$x_{12}$}
 \psfrag{p}[Bl][bl][0.9]{$P$}
 \psfrag{p_1}[Bl][bl][0.9]{$P_1$}
 \psfrag{p_2}[Bl][bl][0.9]{$P_2$}
 \psfrag{p_12}[Bl][bl][0.9]{$P_{12}$}
 \psfrag{LIE}[Bl][bl][0.9]{${\rm LIE}$}
 \psfrag{LAGUERRE}[Bl][bl][0.9]{${\rm LAGUERRE}$}
 \psfrag{MÖBIUS}[Bl][bl][0.9]{${\rm M\ddot{O}BIUS}$}
 \center{\includegraphics[width=130mm]{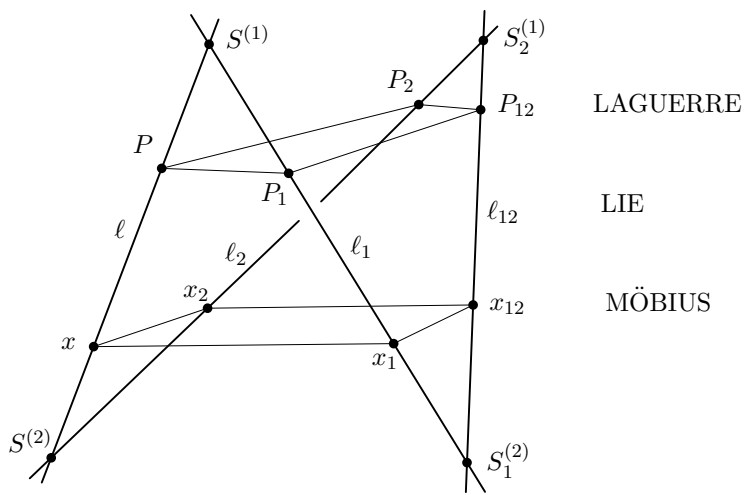}}
 \caption{Elementary quadrilateral of a curvature line parametrized
 surface with vertices $x$ and tangent planes $P$ in the projective model.
 The vertices $x$ build a circular net (M\"obius geometry), and lie in the planes
 $P$ building a conical net (Laguerre geometry). Contact
 elements $(x,P)$ are represented by isotropic lines $\ell$ (Lie
 geometry). Principal curvature spheres $S^{(i)}$ pass through pairs
 of neighboring points $x, x_i$ and are tangent to the
 corresponding pairs of planes $P, P_i$. }
 \label{Fig: Lie quad}
\end{figure}

These relations can be summarized as in Fig.~\ref{Fig: Lie quad}.
Note that the axes of conical nets corresponding to a given
circular net coincide with the Gauss map at its vertices,
considered by Schief \cite{Sch2}.
\smallskip

{\bf Remark.} In the situations of Corollary \ref{cor: synthesis},
i.e., when the R-congruence $S$ consists of points $x$ (and is
therefore a circular net) or of planes $P$ (and is therefore a
conical net), the elementary construction step of Lemma
\ref{lemma: synthesis} allows for a very simple description from
the Euclidean perspective in $\bbR^3$. This has been given by
Pottmann \cite{P}.

\begin{itemize}
\item[i)] {\em Given a contact element $(x,P)$ and a point $x_1$,
find a plane $P_1$ through $x_1$ so that there exists a sphere
$S^{(1)}$ tangent to both planes $P$, $P_1$ at the points $x$,
$x_1$, respectively}. Solution: $P_1$ is obtained from $P$ by the
reflection in the bisecting orthogonal plane of the edge
$[x,x_1]$. The center $c^{(1)}$ of the sphere $S^{(1)}$ is found
as the intersection of the normal to $P$ at $x$ with the bisecting
orthogonal plane of the edge $[x,x_1]$.

\item[ii)] {\em Given a contact element $(x,P)$ and a plane $P_1$,
find a point $x_1$ in $P_1$ so that there exists a sphere
$S^{(1)}$ tangent to both planes $P$, $P_1$ at the points $x$,
$x_1$, respectively}. Solution: the point $x_1$ is obtained from
$x$ by the reflection in the bisecting plane of the dihedral angle
formed by $P$, $P_1$. The center $c^{(1)}$ of the sphere $S^{(1)}$
is found as the intersection of the normal to $P$ at $x$ with this
bisecting plane.
\end{itemize}

\section{R-congruences of spheres} \label{sect: R-congr}

In Sect. \ref{Sect: curv}, Corollaries \ref{cor: curv spheres},
\ref{cor: Ribaucour spheres}, we have seen that principal
curvature spheres of a principal contact element net and spheres
of a Ribaucour transformation build discrete R-congruences,
introduced in Definition \ref{def: R-congr}. In this section we
study the geometry of discrete R-congruences of spheres.
Definition \ref{def: R-congr} can be re-formulated as follows: a
map
\[
S:\bbZ^2\to\{\rm oriented\ spheres\ in\ \bbR^3\},
\]
or the corresponding map
\[
\hat{s}:\bbZ^2\to\bbL^{4,2}\subset\bbR^{4,2}
\]
into the space of homogeneous coordinates, is called a discrete
R-congruence of spheres, if for any $u\in\bbZ^m$ and for any pair
$1\le i\neq j\le m$ the linear subspace
\[
\Sigma={\rm span}(\hat{s},\hat{s}_i,\hat{s}_j,\hat{s}_{ij})
\]
is three-dimensional. Thus, to any elementary square of $\bbZ^m$
there corresponds a three-dimensional linear subspace
$\Sigma\subset\bbR^{4,2}$.

The R-congruence of principal curvature spheres $S^{(i)}$ of the
$i$-th coordinate direction is degenerate in the sense that the
subspaces of its elementary quadrilaterals
\[
\Sigma={\rm
span}(\hat{s}^{(i)},\hat{s}_i^{(i)},\hat{s}_{ij}^{(i)},\hat{s}_j^{(i)})
\]
contain two-dimensional isotropic subspaces (corresponding to
$\ell_i$ and $\ell_{ij}$). The R-congruence of spheres of a
generic Ribaucour transformation is, on the contrary,
non-degenerate: its $\Sigma$'s do not contain two-dimensional
isotropic subspaces, and its elementary quadrilaterals are
included in planar families of spheres, introduced in the
following definition.

\begin{definition}\label{Def: planar family of spheres}
{\bf (Planar family of spheres)} A planar family of spheres is a
set of spheres whose representatives $\hat{s}\in\bbP(\bbL^{4,2})$
are contained in a projective plane $\bbP(\Sigma)$, where $\Sigma$
is a three-dimensional linear subspace of $\bbR^{4,2}$ such that
the restriction of $\langle\cdot,\cdot\rangle$ to $\Sigma$ is
non-degenerate.
\end{definition}
Thus, a planar family of spheres is an intersection
$\bbP(\Sigma\cap\bbL^{4,2})$. Clearly, there are two
possibilities:
\begin{itemize}
\item Signature of $\langle\cdot,\cdot\rangle|_\Sigma$ is $(2,1)$,
so that signature of $\langle\cdot,\cdot\rangle|_{\Sigma^\perp}$
is also $(2,1)$.

\item Signature of $\langle\cdot,\cdot\rangle|_\Sigma$ is $(1,2)$,
so that signature of $\langle\cdot,\cdot\rangle|_{\Sigma^\perp}$
is $(3,0)$.
\end{itemize}
It is easy to see that a planar family is one-parametric,
parametrized by a circle $\bbS^1$. Indeed, if $e_1,e_2,e_3$ is an
orthogonal basis of $\Sigma$ such that $\langle
e_1,e_1\rangle=\langle e_2,e_2\rangle=-\langle e_3,e_3\rangle=1$
(say), then the spheres of the planar family come from the linear
combinations $\hat{s}=\alpha_1e_1+\alpha_2e_2+e_3$ with
\[
\langle
\alpha_1e_1+\alpha_2e_2+e_3,\alpha_1e_1+\alpha_2e_2+e_3\rangle=0
\quad\Leftrightarrow\quad \alpha_1^2+\alpha_2^2=1.
\]
In the second of the cases mentioned above, the space
$\Sigma^\perp$ only has a trivial intersection with $\bbL^{4,2}$,
so that the spheres of the planar family
$\bbP(\bbL^{4,2}\cap\Sigma)$ have no common touching spheres. This
case has no counterpart in the smooth differential geometry. From
the point of view of discrete differential geometry the first case
is more significant.

\begin{definition}
{\bf (Cyclidic family of spheres)} A planar family of spheres is
called cyclidic, if the signature of
$\langle\cdot,\cdot\rangle|_\Sigma$ is $(2,1)$, so that the
signature of $\langle\cdot,\cdot\rangle|_{\Sigma^\perp}$ is also
$(2,1)$.
\end{definition}

Thus, for a cyclidic family $\bbP(\bbL^{4,2}\cap\Sigma)$ there is
a {\em dual} cyclidic family $\bbP(\bbL^{4,2}\cap\Sigma^\perp)$
such that any sphere of the first one is in an oriented contact
with any sphere of the second one. The family
$\bbP(\bbL^{4,2}\cap\Sigma)$, as any one-parameter family of
spheres, envelopes a canal surface in $\bbR^3$, and this surface
is an envelope of the dual family
$\bbP(\bbL^{4,2}\cap\Sigma^\perp)$, as well. Such surfaces are
called {\em Dupin cyclides}. Thus, to any elementary quadrilateral
of a discrete R-congruence whose spheres
$(\hat{s},\hat{s}_i,\hat{s}_{ij},\hat{s}_j)$ span a subspace of
the signature (2,1) there corresponds a Dupin cyclide.

\begin{figure}[htbp]
 \center{\includegraphics[width=60mm]{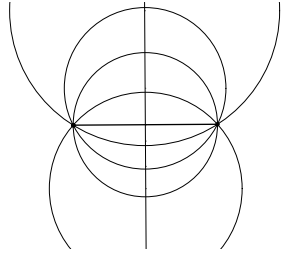}}
 \caption{A cyclidic family of spheres through a circle}
 \label{Fig: circular cyclide}
\end{figure}

Examples:

$\blacktriangleright\ $ points of a circle build a planar cyclidic
family of spheres (of radius zero). The dual family consists of
all (oriented) spheres through this circle, with centers lying on
the line through the center of the circle orthogonal to its plane,
see. Fig.~\ref{Fig: circular cyclide}. The corresponding Dupin
cyclide is the circle itself. It can be shown that any Dupin
cyclide is an image of this case under a Lie sphere
transformation. For a circular net, considered as a discrete
R-congruence, each elementary quadrilateral carries such a
structure.

$\blacktriangleright\ $ planes tangent to a cone of revolution
build a planar cyclidic family of spheres, as well. The dual
family consists of all (oriented) spheres tangent to the cone,
with centers lying on the axis of the cone, see. Fig.~\ref{Fig:
conical cyclide}. The corresponding Dupin cyclide is the cone
itself. For a conical net, considered as a discrete R-congruence,
each elementary quadrilateral carries such a structure.

\begin{figure}[htbp]
 \center{\includegraphics[width=90mm]{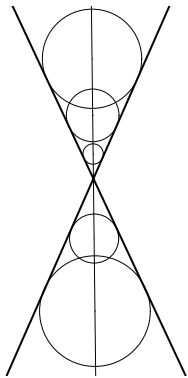}}
 \caption{A cyclidic family of spheres tangent to a cone}
 \label{Fig: conical cyclide}
\end{figure}

\begin{theorem}\label{Thm: two quads of R-congr}
{\bf (Common tangent spheres of two neighboring quadrilaterals of
an R-congruence)} For two neighboring quadrilaterals of a discrete
R-congruence of spheres, carrying cyclidic familes, there are
generically exactly two spheres tangent to all six spheres of the
congruence.
\end{theorem}
{\bf Proof.} Let the quadrilaterals in question belong to the
planar families generated by the subspaces $\Sigma_1$ and
$\Sigma_2$ of the signature (2,1). These quadrilaterals share two
spheres $\hat{s}_1$ and $\hat{s}_{2}$, which span a linear space
of the signature (1,1). Each of the planar families $\Sigma_1$ and
$\Sigma_2$ adds one space-like vector, so that the linear space
$\Sigma_1\cup\Sigma_2$ spanned by all six spheres of the
congruence is four-dimensional and has the signature (3,1), so
that its orthogonal complement $(\Sigma_1\cup\Sigma_2)^\perp$ is
two-dimensional and has the signature (1,1). Intersection of
$\bbL^{4,2}$ with a two-dimensional linear subspace of the
signature (1,1) gives, upon projectivization, exactly two spheres:
indeed, if $e_1,e_2$ form an orthogonal basis of
$(\Sigma_1\cup\Sigma_2)^\perp$ with $\langle
e_1,e_1\rangle=-\langle e_2,e_2\rangle=1$, then the spheres in
this space correspond to $\alpha_1e_1+\alpha_2e_2$ with
\[
\langle \alpha_1e_1+\alpha_2e_2,\alpha_1e_1+\alpha_2e_2\rangle=0
\quad\Leftrightarrow\quad
\alpha_1^2=\alpha_2^2\quad\Leftrightarrow\quad
\alpha_1:\alpha_2=\pm 1. \quad \Box
\]

\begin{figure}[htbp]
 \center{\includegraphics[width=70mm]{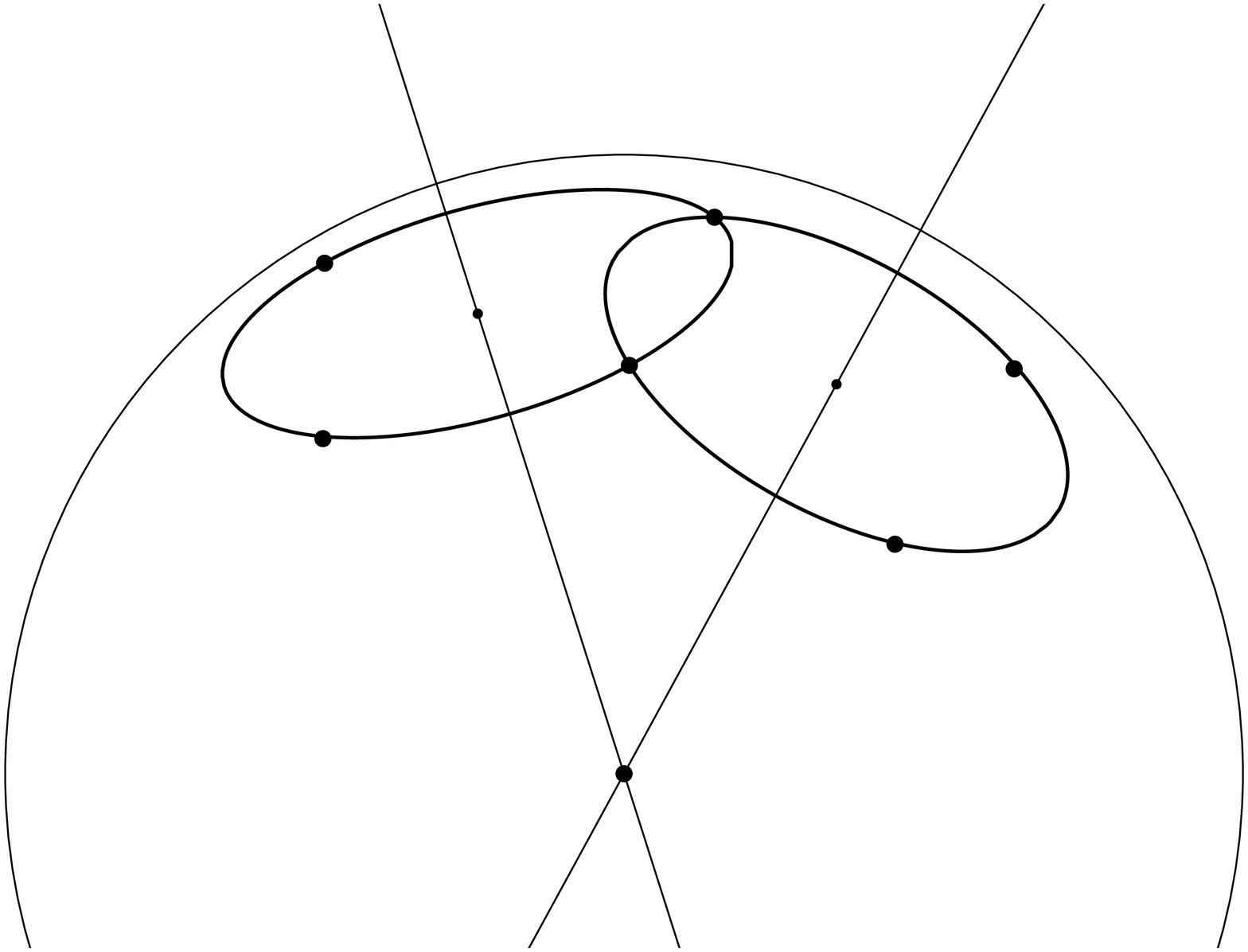}}
 \caption{Normals of two neighboring quadrilaterals of a
  circular net intersect: both lie in the bisecting orthogonal
  plane of the common edge}
 \label{Fig: two circ quads}
\end{figure}

In particular:

$\blacktriangleright\ $ For any two neighboring quadrilaterals of
a circular net, there is one non-oriented sphere (hence two
oriented spheres) containing both circles. Its center is the
intersection point of the lines passing through the centers of the
circles orthogonallly to their respective planes, see
Fig.~\ref{Fig: two circ quads}.

$\blacktriangleright\ $ For any two neighboring quadrilaterals of
a conical net, there is a unique oriented sphere touching both
cones (the second such sphere is the point at infinity). The
center of this sphere is the intersection point of the axes of the
cones.
\begin{figure}[htbp]
 \center{\includegraphics[width=80mm]{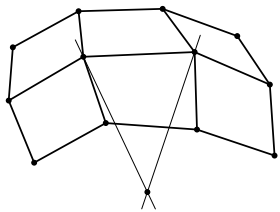}}
 \caption{Normals of two neighboring quadrilaterals of a conical
   net intersect: two common planes of the quadrilaterals are tangent
   to both cones, therefore the axes of both cones
   lie in the bisecting plane of these two planes}
 \label{Fig: two conic quads}
\end{figure}

The next theorem is proved in exactly the same way as Theorem
\ref{Thm: two quads of R-congr}:
\begin{theorem}\label{Thm: hexahedron of R-congr}
{\bf (Common tangent spheres of an elementary hexa\-hedron of an
R-congruence)} For an elementary hexahedron of a discrete
R-congruence of spheres, with all faces carrying cyclidic
families, there are generically exactly two spheres tangent to all
eight spheres at its vertices.
\end{theorem}

It should be mentioned that these spheres, attached to elementary
hexahedra, do not form a discrete R-congruence, contrary to what
has been alleged by A.~Doliwa as a main result of \cite{D6}.
\smallskip

We now turn to a geometric characterization of discrete
R-congruences. From eq. (\ref{eq: Lie sphere}) there follows
immediately that a map
\[
S:\bbZ^m\to\{\rm oriented\ spheres\ in\ \bbR^3\},
\]
is a discrete R-congruence, if and only if the centers
$c:\bbZ^m\to\bbR^3$ form a Q-net in $\bbR^3$, and the two
real-valued functions,
\[
|c|^2-r^2:\bbZ^m\to\bbR\quad {\rm and}\quad r:\bbZ^m\to\bbR,
\]
satisfy the same equation of the type (\ref{eq: dcn property}) as
the centers $c$. By omitting the latter requirement for the signed
radii $r$, one comes to a less restrictive definition than that of
R-congruence. Actually, this definition belongs to M\"obius
geometry and refers to notations of Sect. \ref{subsect: Moebius}
(with $N=3$).

\begin{definition}\label{def: Q-congr}
{\bf (Q-congruence of spheres)} A map
\begin{equation}\label{eq: S-map Moeb}
S:\,\bbZ^m\to\big\{\mbox{\rm non-oriented\ spheres\ in\
}\bbR^3\big\},
\end{equation}
is called an Q-congruence of spheres, if the corresponding map
\begin{equation}\label{eq: s-map Moeb}
\hat{s}:\,\bbZ^m\to\bbP(\bbR^{4,1}_{\rm out}),\quad
\hat{s}=c+\ee_0+\big(|c|^2-r^2\big)\ee_\infty,
\end{equation}
is a Q-net in $\bbP(\bbR^{4,1})$.
\end{definition}
Thus, a map (\ref{eq: S-map Moeb}) is a Q-congruence, if and only
if the centers $c:\bbZ^m\to\bbR^3$ of the spheres $S$ form a Q-net
in $\bbR^3$, and the function $|c|^2-r^2$ satisfies the same
equation (\ref{eq: dcn property}) as the centers $c$.

\begin{theorem}\label{thm: R-congr from Q-congr}
{\bf (Characterization of R- among Q-congruences)} Four (oriented)
spheres $(S,S_i,S_{ij},S_j)$ in $\bbR^3$ comprise an elementary
quadrilateral of an R-congruence, if and only if they comprise (as
non-oriented spheres) an elementary quadrilateral of a
Q-congruence, and satisfy additionally the following condition:
\begin{itemize}
\item[(R)] There exists a non-point sphere in an oriented contact
with all four oriented spheres $S,S_i,S_j, S_{ij}$.
\end{itemize}
Under this condition, any sphere in an oriented contact with three
spheres $S,S_i,S_j$ is in an oriented contact with the fourth one,
$S_{ij}$, as well.
\end{theorem}
{\bf Proof.} Let $S_0$ be a sphere with the center $c_0$ and
(finite) oriented radius $r_0\neq 0$ in an oriented contact with
the three spheres $S, S_i, S_j$. This means that the following
conditions are satisfied:
\begin{equation}\label{eq: tang aux}
\langle c,c_0\rangle
-\tfrac{1}{2}(|c|^2-r^2)-\tfrac{1}{2}(|c_0|^2-r_0^2)-rr_0=0
\end{equation}
(tangency of $S$, $S_0$, cf. (\ref{eq: sph tang})), and two
similar equation with $(c,r)$ replaced by $(c_i,r_i)$ and
$(c_j,r_j)$. Now, using the fact that $c$ and $|c|^2-r^2$ satisfy
one and the same equation of the type (\ref{eq: dcn property}), we
conclude that eq. (\ref{eq: tang aux}) is fulfilled for
$(c_{ij},r_{ij})$, if and only if $r$ satisfies the same equation
(\ref{eq: dcn property}) as $c$ and $|c|^2-r^2$ do. This proves
the theorem in the case when the common tangent sphere $S_0$ for
the three spheres $S, S_i, S_j$ has a finite radius. The case when
$S_0$ has an infinite radius, i.e., is actually a plane, is dealt
with similarly, with the help of equation
\begin{equation}\label{eq: tang aux 1}
\langle c,v_0\rangle-r-d_0=0,
\end{equation}
which comes to replace eq. (\ref{eq: tang aux}). $\Box$
\smallskip

{\bf Remark.} We have already seen that, generically, if three
oriented spheres $S$, $S_i$, $S_j$ have a common sphere in an
oriented contact, then they have a one-parameter (cyclidic) family
of common touching spheres, represented by a three-dimensional
linear subspace $\Sigma$ of $\bbR^{4,2}$. It is easy to see that
if the projection of $\Sigma$ onto $\ee_\infty^\perp$ is
non-vanishing, then the family of spheres represented by
$\Sigma^\perp$ contains exactly two planes. (The only exceptional
case is that of a conical cyclidic family $\Sigma$, all of whose
elements have vanishing $\ee_0$-component and represent planes,
while the family $\Sigma^\perp$ contains no planes.) Therefore, in
all cases but the conical one, condition (R) can be replaced by
the following requirement:
\begin{itemize}
\item[${\rm (R_0)}$] The four oriented spheres $S,S_i,S_j, S_{ij}$
have a common tangent plane (actually, two common tangent planes).
\end{itemize}

It remains to give a geometric characterization of Q-congruences.
This is done in the following theorem.

\begin{theorem}\label{thm: Q-congr}
{\bf (Three types of Q-congruences)} Four (non-oriented)
\linebreak spheres $(S, S_i, S_{ij}, S_j)$ in $\bbR^3$ comprise an
elementary quadrilateral of a Q-congruence, if and only if they
satisfy one of the following three conditions:
\begin{itemize}
\item[(i)] they have a common orthogonal circle, or

\item[(ii)] they intersect along a pair of points (a 0-sphere), or
else

\item[(iii)] they intersect at exactly one point.
\end{itemize}
Case (iii) can be regarded as a degenerate case of both (i) and
(ii).
\end{theorem}

\noindent {\bf Conceptual proof.} {\em Caution:} notations in this
proof refer to the M\"obius-geometric objects which are different
from the Lie-geometric objects denoted by the same symbols. The
linear subspace $\Sigma$ of $\bbR^{4,1}$ spanned by the points
$\hat{s}$, $\hat{s}_i$, $\hat{s}_j$, $\hat{s}_{ij}$ is
three-dimensional, so that its orthogonal complement
$\Sigma^\perp$ is two-dimensional. If $\Sigma^\perp$ lies in
$\bbR^{4,1}_{\rm out}$, i.e., if the restriction of the Minkowski
scalar product to $\Sigma^\perp$ is positive-definite (of
signature (2,0)), then $\Sigma^\perp$ represents a 1-sphere (a
circle) orthogonal to our four spheres, and we have the case (i).
If, on the contrary, the restriction of the scalar product to
$\Sigma^\perp$ has signature $(1,1)$, so that $\Sigma$ lies in
$\bbR^{4,1}_{\rm out}$, then $\Sigma$ represents a $0$-sphere
which is the intersection of our four spheres, and we have the
case (ii). Finally, if the restriction of the scalar product to
$\Sigma$ is degenerate, then $\Sigma\cap\Sigma^\perp$ is an
isotropic one-dimensional linear subspace, which represents the
common point of our four spheres, and we have the case (iii).
$\Box$
\smallskip

\begin{figure}[htbp]
 \psfrag{c1}[Bl][bl][0.9]{$c_j$}
 \psfrag{c2}[Bl][bl][0.9]{$c_{ij}$}
 \psfrag{c3}[Bl][bl][0.9]{$c_i$}
 \psfrag{c4}[Bl][bl][0.9]{$c$}
 \psfrag{r1}[Bl][bl][0.9]{$r_j$}
 \psfrag{r2}[Bl][bl][0.9]{$r_{ij}$}
 \psfrag{r3}[Bl][bl][0.9]{$r_i$}
 \psfrag{r4}[Bl][bl][0.9]{$r$}
 \psfrag{C}[Bl][bl][0.9]{$C$}
 \center{\includegraphics[height=70mm]{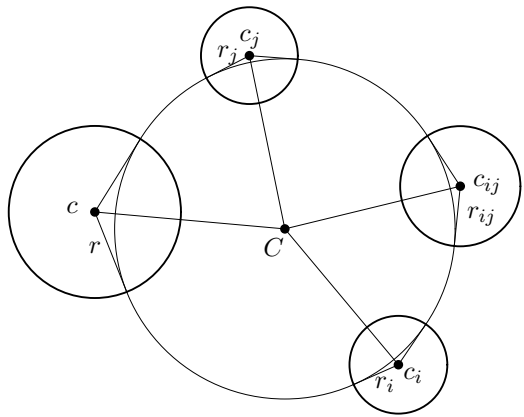}}
 \caption{Elementary quadrilateral of a Q-congruence of spheres,
 the orthogonal circle case}
 \label{Fig: S-circ}
\end{figure}

\noindent {\bf Computational proof.} The quadrilateral in $\bbR^3$
with the vertices at the sphere centers $c$, $c_i$, $c_j$,
$c_{ij}$ is planar; denote its plane by $\Pi$. In the same way as
in the proof of Theorem \ref{thm: dos} we show that there is a
point $C\in\Pi$ such that
\begin{equation}\label{eq: S-ortho}
|c-C|^2-r^2=|c_i-C|^2-r_i^2=|c_j-C|^2-r_j^2=|c_{ij}-C|^2-r_{ij}^2.
\end{equation}
Indeed, the first two of these equations define $C$ uniquely as
the intersection of two lines $\ell_i$ and $\ell_j$ in $\Pi$,
where
\[
\ell_i=\{x\in\Pi: \langle 2x-c_i-c,c_i-c\rangle=r^2-r_i^2\},
\]
and then the last equation in (\ref{eq: S-ortho}) is automatically
satisfied. If the common value of all four expressions in
(\ref{eq: S-ortho}) is positive (say, equal to $R^2$), then the
four spheres under consideration are orthogonal to the circle in
the plane $\Pi$ with the center $C$ and radius $R$, so that we
have the case (i), see Fig. \ref{Fig: S-circ}. If the common value
of (\ref{eq: S-ortho}) is negative (say, equal to $-R^2$), then
the pair of points on the line through $C$ orthogonal to $\Pi$, at
the distance $R$ from $C$, belong to all four spheres, so that we
have the case (ii). Finally, if the common value of (\ref{eq:
S-ortho}) is equal to 0, then $C$ is the intersection point of all
four spheres, and we have the case (iii). $\Box$.
\medskip

Clearly, case (i) of Q-congruences reduces to circular nets, if
the radii of all spheres become infinitely small, cf.
Fig.~\ref{Fig: S-circ}. Q-congruences with intersections of type
(ii) are natural discrete analogs of sphere congruences
parametrized along principal lines. \smallskip

Some remarks about Q-congruences of spheres are in order here.
They are multidimensionally consistent, with the following
reservation: given seven points $\hat{s}$, $\hat{s}_i$,
$\hat{s}_{ij}$ in $\bbP(\bbR^{4,1}_{\,\rm out})$, the Q-property
(planarity condition) uniquely defines the eighth point
$\hat{s}_{123}$ in $\bbP(\bbR^{4,1})$, which, however, might get
outside of $\bbP(\bbR^{4,1}_{\,\rm out})$, and therefore might not
represent a real sphere. Thus, the corresponding discrete 3D
system is well-defined on an open subset of the space of initial
data only. As long as it is defined, it can be used to produce
transformations of Q-congruences, with usual permutability
properties.

Note the following difference between Q-congruences and
R-congruences: given three spheres $S$, $S_i$, $S_j$ of an
elementary quadrilateral, one has a two-parameter family for the
fourth sphere $S_{ij}$ in the case of a Q-congruence, and only a
one-parameter family in the case of an R-congruence. This is a
consequence of the fact that $\bbR^{4,1}_{\rm out}$ is an open set
in $\bbR^{4,1}$, while $\bbL^{4,2}$ is a hypersurface in
$\bbR^{4,2}$.

\begin{appendix}

\section{Appendix: cyclographic model of Laguerre \\ geometry}
\label{Subsect: Laguerre cyclographic}

In the cyclographic model of Laguerre geometry, the preferred
space is the space of hyperspheres $(\bbR^{N,1,1})^*$, so
hyperspheres $S\subset\bbR^N$ are modelled as points
$\hat{s}\in\bbP\big((\bbR^{N,1,1})^*\big)$, while hyperplanes
$P\subset\bbR^N$ are modelled as hyperplanes $\{\xi: \langle
\hat{p},\xi\rangle=0\}\subset\bbP\big((\bbR^{N,1,1})^*\big)$.
Thus, a hyperplane $P$ is interpreted as a set of hyperspheres $S$
which are in oriented contact with $P$.

Basic features of this model:
\begin{itemize}
\item[(i)] The set of oriented hyperspheres $S\subset\bbR^N$ is in
a one-to-one correspondence with points
\begin{equation}\label{eq: Lag sphere cyclogr}
\sigma=(c,r)
\end{equation}
of the Minkowski space $\bbR^{N,1}$ spanned by the vectors
$\ee_1,\ldots,\ee_N,\ee_{N+3}$. This space has an interpretation
of an affine part of $\bbP\big((\bbR^{N,1,1})^*\big)$.

\item[(ii)] Oriented hyperplanes $P\subset\bbR^N$ can be modelled
as hyperplanes in $\bbR^{N,1}$:
\begin{equation}\label{eq: Lag plane cyclogr}
\pi=\big\{(c,r)\in\bbR^{N,1}:\,\langle (v,1),(c,r)\rangle=\langle
v,c\rangle-r=d\big\}.
\end{equation}
Thus, oriented hyperplanes $P\in\bbR^N$ are in a one-to-one
correspondence with hyperplanes $\pi\subset\bbR^{N,1}$ which make
angle $\pi/4$ with the subspace
$\bbR^N=\{(x,0)\}\subset\bbR^{N,1}$.

\item[(iii)] An oriented hypersphere $S\subset\bbR^N$ is in an
oriented contact with an oriented hyperplane $P\subset\bbR^N$, if
and only if $\sigma\in\pi$.

\item[(iv)] Two oriented hyperspheres $S_1,S_2\subset\bbR^N$ are
in an oriented contact, if and only if their representatives in
the Minkowski space $\sigma_1,\sigma_2\in\bbR^{N,1}$ differ by an
isotropic vector: $|\sigma_1-\sigma_2|=0$.
\end{itemize}
In the cyclographic model, the group of Laguerre transformations
admits a beautiful description:
\begin{theorem}\label{thm: Lag fund}
{\bf (Fundamental theorem of the Laguerre geometry)} The group of
Laguerre transformations is isomorphic to the group of affine
transformations of $\bbR^{N,1}$: $y\mapsto \lambda Ay+b$, where
$A\in O(N,1)$, $\lambda>0$, and $b\in\bbR^{N,1}$.
\end{theorem}

\end{appendix}

\end{document}